\documentclass[10pt]{article}
\usepackage[a4paper]{geometry}
\usepackage{url}
\usepackage{fullpage}
\usepackage{amsmath}
\usepackage{amsfonts}
\usepackage{amssymb}
\usepackage{hyperref}
\usepackage{times}

 \newcommand{\balpha}{\boldsymbol{\alpha}}
 \newcommand{\btheta}{\boldsymbol{\omega}}
  \newcommand{\be}{\boldsymbol{e}}
  \newcommand{\bx}{\boldsymbol{x}}
\newcommand{\by}{\boldsymbol{y}}
\newcommand{\bz}{\boldsymbol{z}}
  \newcommand{\bu}{\boldsymbol{u}}
  \newcommand{\bv}{\boldsymbol{v}}

\newcommand{\cS}{{\cal{S}}}

\newcommand{\cC}{{\mathcal{C}}}

\newcommand{\cB}{{\mathcal{B}}}

\newcommand{\eop}{\hfill $\Box~~$\\[1mm]}
\newcommand{\Pf}{{\em Proof:~~}}

\newcommand{\bi}{\begin{itemize}} 
\newcommand{\ei}{\end{itemize}}
\newcommand{\beq}{\begin{equation}} 
\newcommand{\ee}{\end{equation}}

\usepackage{fancyhdr}
\usepackage{url}
\usepackage{amsmath}
\usepackage{amsthm}
\usepackage{amsfonts}
\usepackage[T1]{fontenc}
\usepackage[latin1]{inputenc}
\usepackage{amssymb}%
\usepackage{graphicx}
\usepackage{amsmath} 
\newtheorem{definition}{Definition}
\newtheorem{example}{Example}

\newtheorem{lemma}{Lemma}

\newtheorem{theorem}{Theorem}
\newtheorem{proposition}{Proposition}

\newcommand{\cprank}[1]{{ \emph{CP-rank} }}

\usepackage{fancyhdr}
\fancypagestyle{plain}{ %
\fancyhf{} 
}

\title{Approximation Hierarchies for Copositive Tensor Cone}
\author{
Muhammad Faisal Iqbal $^{1,2}$
\thanks{Corresponding author: 
1. Allama Iqbal Open University, Islamabad Pakistan, Department of Mathematics, ~ Email: \url{faisal.iqbal@aiou.edu.pk}  ~
2. Institute of Space Technology, Islamabad, Pakistan, Department of Applied Mathematics, ~ Email: \url{fasi_g88@yahoo.com} } \and
Faizan Ahmed\thanks{\small Institute of Space Technology, Islamabad, Pakistan, Department of Applied Mathematics,   Email: \url{faizan.ahmed@mail.ist.edu.pk}} \and Muhammad Aqeel \thanks{\small Institute of Space Technology, Islamabad, Pakistan, Department of Applied Mathematics,   Email: \url{muhammad.aqeel@mail.ist.edu.pk}} \and Salman Ahmad\thanks{\small Institute of Space Technology, Islamabad, Pakistan, Department of Applied Mathematics,   Email: \url{salman.ahmad@mail.ist.edu.pk}} 
}

\begin{document}
\markboth{Author}{Title of the paper}
\date{\today}
\maketitle
\begin{abstract} 
\noindent
In this paper we discuss copositive tensors, which are a natural generalization of the copositive matrices. We present an analysis of some basic properties of copositive tensors; as well as the conditions under which class of copositive tensors and the class of positive semidefinite tensors coincides. Moreover, we have describe several hierarchies that approximates the cone of copositive tensors. The hierarchies are predominantly based on different regimes such as; simplicial partition, rational griding and polynomial conditions. The hierarchies approximates the copositive cone either from inside (inner approximation) or from outside (outer approximation). We will also discuss relationship among different hierarchies.\\  

\end{abstract}
\section{Introduction} 
Copositive optimization has become an active area of research during recent years. The significance of copositive programming is due to the fact that; several combinatorial and non-convex optimization problems can have linear programming reformulation over the copositive cone, which is convex. A growing list of problems that have copositive programming reformulations includes; standard quadratic programming (\cite {bomze2000copositive}, \cite{bome01solving}), the chromatic and stability number of a graph (\cite{Klerk02app}, \cite{JVPena07comp},\cite{Monique08theOpertr} \cite{DIFranz10copoPro}), crossing number of a graph \cite{de2006improved}, the maximum stable set problem \cite{de2006improved}, the quadratic assignment problem \cite{Povh:2009:CSR:2296586.2296766} and discrete optimization\cite{burer2009copositive}. Consequently the new developments about copositive cone and its dual cone can be helpful to the solution of all the above mentioned hard problems. \\
The copositive programming problems are not solvable directly, as the copositive and completely positive cones are not tractable, thus the approximation hierarchies for these cones have been studied in much detail by the matrix theorists, see \cite{Laurent2009}. Several approximation hierarchies based on sum-of-squares conditions and discretization methods have been studied. For instance, Parrilo \cite{Parrilo00structuredsemidefinite} had provided a hierarchy of linear and semi-definite inner approximations for copositive cone (see also\cite{deklerk2002}). Moreover, Bomze and de Klerk \cite{bomze2002solving} (see also \cite{Bundfuss2009}) suggested a criteria to check membership of a given matrix in copositive cone by using a sequence of polyhedral approximations; which approximates the cone from inside and outside. Moreover, these approximations are exact in the limit. 
Problems stated above have a common feature that they have quadratic objective function or in certain cases quadratic constraints. The immediate generalization of quadratic optimization is  polynomial optimization. In recent years, polynomial optimization has attracted many researchers due to its vast applications in empirical modeling of science \cite{Lasserre2012polyOpt} and engineering problems such as biomedical engineering \cite{barmpoutis:17,ghosh:hal-00340600,Zhang2011}, signal processing \cite{Meng:2009:QPA:1653465.1653481,Qi2003,Weiland:2010:SVD:1771983.1772002}, quantum graphs \cite{quantumGraph2015}, and material science \cite{Soare2008915}. 
Quadratic functions can be represented using matrices, likewise polynomial can be represented by multidimensional arrays known as tensors. Similar to the case of quadratic optimization a reformulation of polynomial optimization is described by~\cite{Lasserre2012polyOpt,Pena2015}. The notion of copositive and completely positive tensors is used to describe these reformulations. The copositive and completely positive cones of matrices, which are tensors of order two, are very well explored, therefore it seems natural to study analogous results for copositive tensors. However this generalization is not trivial since higher dimension usually destroy the nice structure present at lower dimension. \\ 
The area is not very well explored, however there are some research in describing the properties of tensors. Song have provided a characterization of copositive tensors using eigenvectors of principal sub-tensors~\cite{song2015necessary}. Qi extended the diagonal dominance sufficient conditions for complete positivity to the tensor case~\cite{qi2014nonnegative}. In this article, we analyze some approximation schemes for copositive cone of tensors. Therefore in order to answer the theoretical questions posed above, the contributions of this article are: (a) The basic properties of copositive tensor cone $\mathcal{C}_{n,d}$ are discussed. Moreover, the conditions, for which copositive cone of tensors coincides with the positive semidefinite cone of tensors, are established.
(b) Several polynomial conditions based approximation hierarchies for copositive cone of tensors are presented, which are generalization of analogous results for matrix case.
(c) Approximation hierarchies for copositive cone of tensors; which are based on simplicial partition and rational griding are also presented. (d) It has been established that the above mentioned hierarchies approximates the cone $\mathcal{C}_{n,d}$ exactly, in the limiting case. (e)The inclusion relations among approximation hierarchies are also presented.

The article is arranged as; Section 2 comprises of the basic definitions and notations. In Section 3, first we give a brief introduction of tensor cones, secondly we discuss several properties of tensors along with characterization of tensor cones. In Section 4.1, we present the inner approximation hierarchies for copositive cone of tensors. One of these hierarchies is based on polynomial conditions and the other is based on simplicial partition. Then we present the containment relationship among inner approximation hierarchies for copositive cone of tensors. Lastly in Section 4.2, two types of outer approximation hierarchies for copositive cone of tensors are presented based on polynomial conditions and rational griding. In section 5 we provide conclusion and future work.
\section{Preliminaries} 
The use of matrix for the representation of quadratic form has become ubiquitous, in recent years. Thus it looks natural to represent a polynomial using multi-dimensional array, usually termed as tensor. Throughout this article; $\Re^n$ denotes the $n$-dimensional Euclidean space and $\Re^n_+$ denotes the non-negative orthant of $\Re^n$. The set of natural numbers is denoted by $\mathbb{N}$ and the set of whole numbers is denoted by $\mathbb{N}_{0}=\{0,1,2,\cdots\}$. The vectors are denoted by using small case bold letters and matrices are denoted by using the capital letters, however calligraphic capital letters are used to denote tensors.  The tensor is defined as follows;
\begin{definition}[Tensor]
	A Tensor is a multi-array of real numbers. Mathematically, a Tensor; \[\mathcal{A}=\big(a_{{i_{1}} \ldots {i_{d}} } \big)_{1\leq i_{1}, \ldots, i_{d} \leq n}.\]  is a $n$-dimensional, $d^{th}$-order array; and $\mathcal{A}$ is said to be symmetric if all permutations $\sigma$ of  indices's represents the same element of $\mathcal{A}$,\ i.e we have
	\[a_ {i_1 \ldots i_d} = a_ {\sigma({i_{1} \ldots i_{d}})}\] 
\end{definition}
\noindent
Clearly any matrix is a tensor of order two. For brevity of notation, if some index $i_j$ of an element $a_{i_1i_2 \cdots i_d} \in \mathcal{A}$ is repeated $k$-times, we write it as ${(i_j)}^{k}$ and such elements are denoted by \[a_{{\underset{k-times}{\underbrace{i_ji_j \cdots i_j}}}i_{(k+1)}i_{(k+2)}\cdots i_{d}}=a_{{(i_j)^k}i_{(k+1)}i_{(k+2)} \cdots i_{d}}\]
The collection of $n$-dimensional, $d^{th}$-order symmetric is denoted by $\mathcal{S}_{n,d}$. One particular case is the cone of entry-wise non-negative tensors denoted by $\mathcal{N}_{n,d}$. For $ \bx \in \Re^{n}$ the product tensor $\mathcal{X}={\bx}^d$; is stated as under :
\begin{align*}
\mathcal{X}=\underset{d-times}{\underbrace{{\bx} \otimes \cdots \otimes {\bx}}} \ \in  {({\Re}^n \otimes \cdots \otimes {\Re}^n)}
\end{align*} Let $\mathcal{T}_d:\Re^n\rightarrow\mathcal{S}_{n,d}$ be a mapping defined as; $\mathcal{T}_d(\bx)=\mathcal{X}$.
\\
Using above notation the $d^{th}$ degree homogeneous polynomial in $n$-variables is stated as,  
\begin{align}\label{2}
f_\mathcal{A}{({\bx})} =  \sum_{i_1,i_2 \cdots, i_d =1}^{n} a_{i_1 \cdots i_d} {x_{i_1} x_{i_2}} \cdots {x_{i_d}}
\end{align}
where $\mathcal{A}$ is $n$-dimensional, $d^{th}$-order symmetric tensors. Introducing a more convenient notation for homogeneous polynomial $f_\mathcal{A}{({\bx})}$ associated with tensor $\mathcal{A}$, denoting  $\mathbb{Z}_d=\{0,1,2,\cdots,d\}$ and for any vector ${\balpha}\in\mathbb{Z}_d^n$, let us define 1-norm $\|{\balpha}\|_{1}=\sum_{i=1}^{n} \alpha_i$. For the subset, $\mathbb{I}^n(d)=\{{\balpha}\in\mathbb{Z}_d^n:\|{\balpha}\|_{1}=d \}$   of $\mathbb{Z}_{d}^n$, the monomial of degree  $d$ over $\Re^n$ is defined as; 
\begin{align}\label{3}
{\bx}^{\balpha} &= \prod_{i=1}^{n} {x_i}^{\alpha_i} \\ &={x_1}^{\alpha_1}{x_2}^{\alpha_2}\cdots{x_n}^{\alpha_n} ~\text{for}~\balpha \in \mathbb{I}^n(d)~\text{and}~ \bx \in \Re^n
\end{align} 
The collection of all possible $d^{th}$-degree monomials in $n$ variables is denoted by $S(\mathcal{X})$, and is given as under:  
\begin{align*}
S(\mathcal{X})=\bigg\lbrace {\bx}^{\balpha} : \balpha \in \mathbb{I}^n(d) ~and~\bx \in \Re^n \bigg\rbrace 
\end{align*}
 Since for an arbitrary $\bx \in \Re^n$ there exists a bijection $\phi_{\bx}: \mathbb{I}^n(d) \rightarrow S(\mathcal{X})$ such that; $\phi_{\bx}(\balpha)=\bx^{\balpha}$.  Thus it implies that, the cardinality of set $S(\mathcal{X})$ is same as the cardinality of $\mathbb{I}^n(d)$, that is $|S(\mathcal{X})|=|\mathbb{I}^n(d)|= {n+d-1 \choose d}$ \cite{Monique2013}. \\
Note that, the set $S(\mathcal{X})$ is the collection of all possibly distinct elements of $\mathcal{X}$. The inner product of vectors $\bu,\bv \in \Re^{n} $ denoted by $ \langle \bu,\bv \rangle $ is defined as; $\langle \bu,\bv \rangle=\bu^{T}\bv$. Moreover, for tensors $\mathcal{A}, {\cB} \in {\cS}_{n,d}$ the inner product $\big\langle \mathcal{A}, {\cB} \big\rangle$ is defined as;
\begin{align}\label{4}
\big\langle \mathcal{A} , \mathcal{B} \big\rangle = \sum_{i_1,i_2 \cdots, i_d =1}^n a_{i_1i_2 \cdots i_d}  b_{i_1i_2 \cdots i_d}
\end{align}
\noindent
 by using (\ref{3}) and (\ref{4}), we may rewrite (\ref{2}) as; 
\begin{align*}
f_\mathcal{A}{({\bx})} &=\sum_{i_1,i_2 \cdots, i_d =1}^{n} a_{i_1i_2 \cdots i_d} x_{i_1} x_{i_2} \cdots x_{i_d}\\
&=\sum_{i_1,i_2 \cdots, i_d =1}^{n} a_{i_1i_2 \cdots i_d}\prod_{k \in \mathbb{Z}_d \backslash \{0\}}^{} \bx^{\be_{i_k}} = \bigg\langle \mathcal{A}, {\mathcal{T}_{d}(\bx)} \bigg\rangle
\end{align*}
\begin{align*}
\text{where} \ \be_{i_k}\in \Re^n \text{ is a unit vector having all its components} \ 0 \ \text{except the} \ {i_k}^{th} \ \text{component.} 
\end{align*}
 Let $V$ be a vector space with underline field $F$, for arbitrary vectors $\bu_{1},\bu_{2},\cdots,\bu_{m} \in V$ and scalars $\lambda_{1}, \lambda_{2},\cdots,\lambda_{m}\in F$ the linear combination $\bu=\sum_{i=1}^{m}\lambda_{i}\bu_{i}$ is said to be; affine combination if $\sum_{i=1}^{m}\lambda_i=1$, and it is said to be conical combination if  $\lambda_i \ge 0$. Moreover, $\bu$ is said to be convex combination if it is both affine and conical combination.   The set $\cC \subseteq V$ is said to be a cone if for each $\bx\in \cC$ it implies that $\lambda \bx \in \mathcal{C}$ for all scalars $\lambda \ge 0$. Moreover, the cone $\cC$ is said to be convex cone if for each pair $\bx,\by\in \cC$ and for non-negative scalars $\lambda_1,\lambda_2 \in F$ we have, ${\lambda_1} \bx+{\lambda_2} \by \in \mathcal{C}$. The dual of cone $\mathcal{C}$ denoted by $\cC^*$ is stated as under : 
\begin{align*}
\cC^*=\bigg\lbrace \bu \in V:\big\langle \bu,\bv \big\rangle \ge 0 \; \forall\;  \bu \in \cC \bigg\rbrace
\end{align*} 
For any subset ${\mathcal{M}} \subseteq V$ the conic hull of $\mathcal{M}$ denoted by $conic(\mathcal{M})$ is defined as: 
\[conic(\mathcal{M})=\bigg\lbrace \sum_{i=1}^{m}\lambda_{i}\bu_{i}: \ \bu_{i}\in\mathcal{M},  \ \lambda_{i} \ge 0 \ \ \forall i=1,2,\cdots,m \bigg\rbrace \]
A convex cone $\cC$ is said to be pointed if $\{\cC\}\cap\{-\cC\}=\{0\}$, and $\cC$ is said to be solid if its interior is nonempty. A convex cone which is closed, pointed and solid is termed as proper cone. 
\noindent
A convex cone $\mathcal{C}$ is said to be a polyhedral cone if it is finitely generated, that is, there exists a finite set $\mathcal{M}$ such that, $\mathcal{C}=conic(\mathcal{M}).$ 
\section{Tensor Cones}
In this section, we define several cones of tensors. These cones appears as a generalization of the cone of matrices. We discuss various properties of these cones together with special cases where these cones coincides. The collection of symmetric tensors $\mathcal{S}_{n,d}$ is a vector space over the field of reals $\Re$.  we define our first cone of tensors as follows; 
\begin{definition}[Positive Semidefinite Cone: $\mathcal{S}_{n,d}^+$ ]
	A tensor $\mathcal{A}\in {\cS}_{n,d}$  is said to be positive semidefinite (PSD) if $f_\mathcal{A}(\bx)=\big\langle \mathcal{A},\mathcal{T}_d(\bx) \big\rangle \ge 0$ for $\bx\in \Re^n$. The set of $n$-dimensional, $d^{th}$-order positive semidefinite symmetric tensors denoted by $\mathcal{S}_{n,d}^+ $ is stated as under : 
	\begin{align}\label{PSD}
	\mathcal{S}_{n,d}^{+}   & := \bigg\lbrace \mathcal{A} \in \mathcal{S}_{n,d} : \big\langle \mathcal{A},\mathcal{T}_d(\bx) \big\rangle \ge 0 ~\forall~ \bx \in \Re^n  \bigg\rbrace 
	\end{align}
\end{definition}
\noindent
Clearly $\mathcal{S}_{n,d}^{+}$ is a convex cone, as for each tensor $\mathcal{A}\in \mathcal{S}_{n,d}^{+}$ it implies that $\lambda\mathcal{A}\in\mathcal{S}_{n,d}^{+}$ for $\lambda\in\Re_+$, moreover for $\mathcal{A},\mathcal{B}\in\mathcal{S}_{n,d}^{+}$ we have ${\lambda_1}\mathcal{A}+{\lambda_2}\mathcal{B}\in \mathcal{S}_{n,d}^{+}$ for ${\lambda_1},{\lambda_2}\in\Re_+$. \\
It is clear from the above definition that $d$ must be even, since for odd degree tensors non-negativity requirement can not be satisfied. 
The dual of $\mathcal{S}^{+}_{n,d}$ is usually termed as completely positive semidefinite cone denoted by $ \mathcal{S}^{{+}^{*}}_{n,d} $ is stated as: 
\begin{align*}
\mathcal{S}^{{+}^{*}}_{n,d}&=\bigg\lbrace \mathcal{X} \in \mathcal{S}_{n,d} : \big\langle \mathcal{A},\mathcal{X} \big\rangle \ge 0 ~\forall~ \mathcal{A} \in \mathcal{S}^{+}_{n,d}   \bigg\rbrace\\
&=\bigg\lbrace \sum_{\mathcal{X} \in \mathcal{S}_{n,d}}^{} \mathcal{X}
: \mathcal{X}=\mathcal{T}_{d}(\bx), ~~\forall~\bx\in \Re^n  \bigg\rbrace
\end{align*}
Clearly for $d=2$,  PSD cone is self dual, that is, $\mathcal{S}^{+}_{n,2} = \mathcal{S}^{{+}^{*}}_{n,2}$. However for $d \ge 4$ self duality property is not true, that is $\mathcal{S}^{+}_{n,d} \ne \mathcal{S}^{{+}^{*}}_{n,d}$ in general (see [counter example 4.5,\cite{Luo2015}]).\\
 The polynomials are continuously differentiable functions. Since we know that the Hessian matrix of a function is PSD if and only if it is a convex function[Theorem 4.5. \cite{rockafellar1970convex}]. Therefore the convexity of homogeneous polynomial defined in (\ref{2})  amounts to check if $\triangledown^2f_\mathcal{A}(\bx)\in\mathcal{S}^+_{n,2}~for~\bx\in \Re^n$. In [Proposition 5.10, \cite{Luo2015}] it has been shown that, if the polynomial $f_\mathcal{A}(\bx)$ is convex then its associated tensor $\mathcal{A}$ is PSD, but the converse is not true in general, (see [counter example 5.11, \cite{Luo2015}] and \cite{Ahmadi2012}).\\
 During recent years, copositive matrices has been studied extensively due to their usefulness to solve combinatorial and quadratic optimization problems \cite{Mirjam2010survey} as discussed in the introduction. We define positive semi-definiteness of tensor $\mathcal{A}$ over the non-negative restriction(subset) of $\Re^n$ as follows. 
\begin{definition} [Copositive Cone: ${\cC}_{n,d}$]
	A tensor $\mathcal{A}\in {\cS}_{n,d}$ is said to be copositive if $\big\langle \mathcal{A},\mathcal{T}_d(\bx) \big\rangle \ge 0$ for all $\bx\in \Re^n_{+}$. The set of $n$-dimensional, $d^{th}$-order copositive tensors denoted by $\mathcal{C}_{n,d}$ is stated as under : 
	\begin{align}
	\mathcal{C}_{n,d}   := \bigg\lbrace \mathcal{A} \in {\cS}_{n,d} : \big\langle \mathcal{A},\mathcal{T}_d(\bx) \big\rangle \ge 0 ~\forall~ \bx \in \Re^n_+  \bigg\rbrace
	\end{align}	
\end{definition}
\noindent
Similar to (\ref{PSD}), it is obvious that, $\mathcal{C}_{n,d}$ is also convex cone.
Here we would like to remark that in case of copositive tensor one need not to take $d$ to be even. It is obvious from the above definitions, if a tensor is positive semidefinite then it is copositive as well; however,  copositive tensors are not necessarily positive semidefinite  in general; we present a counter example as under.
\begin{example}
	For a symmetric tensor $\mathcal{A}\in\mathcal{S}_{3,4}$ with entries;
	\begin{align*}
	a_{i_1i_2i_3i_4}= 
	\begin{cases}
	\begin{array}{cc}
	0 & ~~ ~~if~i_{j}=1~\forall~j \in \{1,2,3,4\} \\
	1 & ~~ ~~if~i_{j}=2~\forall~j \in \{1,2,3,4\} \\
	1 & ~~ ~~if~i_{j}=3~\forall~j \in \{1,2,3,4\} \\
	5 & otherwise
	\end{array}
	\end{cases}
	\end{align*}
	The associated polynomial is: 
	\begin{align*}
	f_{\mathcal{A}}(\bx)=\bigg(x_2^4+x_3^4+\sum_{i_1,i_2,i_3,i_4=1}^{3}5x_{i_1}x_{i_2}x_{i_3}x_{i_4}\bigg) ~\ge 0~~\forall~\bx\in\Re^n_+
	\end{align*}
	for $\bar{\bx}=(-2,0,1)^{T}$ we have, $f_{\mathcal{A}}(\bar{\bx})=-24$, which implies that, $\mathcal{A}$ is not positive semidefinite tensor.
\end{example} 
\noindent
We describe some basic properties of copositive tensor $\mathcal{A} \in \mathcal{C}_{n,d}$, which are generalizations of the similar properties for the matrix case.  
\begin{proposition}
	Let $\mathcal{A} \in \mathcal{S}_{n,d}$ be any copositive tensor. Then the following properties holds;
	\begin{itemize}
		\item[(i)] $a_{(i)^d} \ge 0$ for all $i$.
		\item[(ii)] If $a_{(i_j)^d} =0 $ then $a_{i_1i_2 \cdots i_{j}i_{j+1} \cdots i_{d}} \ge 0$ where $i_j \ne i_k$ for all $j \ne k \in \{1,2, \cdots, d \}$
	\end{itemize}
\end{proposition}
\noindent
\Pf For an arbitrary tensor $\mathcal{A} \in \mathcal{C}_{n,d}$.
\begin{itemize}
	\item[(i)] Let ${\be}_i \in \Re^n_+$ be a standard unit vector. Then the tensor $\mathcal{B}=\mathcal{T}_{d}(\be_i)$  has zero elements everywhere except the $i^{th}$ diagonal element i.e. $b_{(i)^d}=1$. Since $\mathcal{A}$ is copositive thus, the inner product $0 \le \big\langle \mathcal{A},\mathcal{T}_{d}({\be}_i) \big\rangle =a_{(i)^d}$; that is, $a_{(i)^d}\ge 0~~\forall~i$.
	\item[(ii)]  let, $\mathcal{A}\in \mathcal{C}_{n,d}$ be a tensor with $a_{{(i_{1})}^d}=0$, we assume on the contrary that there exist ${i_j} \in \{1,2,\cdots,n\}$ such that; 
	\[ a_{{(i_{1})}^{(d-m)}\cdot {(i_{j})}^{(m)}} < 0 \ where \ {i_j} \ne {i_1} \ and \ 1\le m < d \] 
	Then for $\bar{\bx} \in \Re^n_+$ with components; 
	\begin{align*}
	x_{i}=\begin{cases}
	\begin{matrix}
	1 &~~ if ~~i\in \{1,j\} \\  
	0 &~~ otherwise
	\end{matrix}
	\end{cases} 
	\end{align*}
	 the associated homogeneous polynomial $f_\mathcal{A}(\bar{\bx})$ is defined as: 
	\begin{align}\label{jthSum}
	f_{\mathcal{A}}(\bar{\bx})= a_{{(i_{1})}^{(d-m)}\cdot {(i_{j})}^{(m)}}\cdot(1)^{(d-m)}\cdot(1)^{(m)} <0
	\end{align}
	from (\ref{jthSum}), we have a contradiction, i.e. $\mathcal{A}$ is not copositive.
\end{itemize} 
 \eop
The dual of copositive cone $\mathcal{C}_{n,d}$ is completely positive cone; denoted by $\mathcal{C}_{n,d}^*$ and is defined below.  
\begin{definition}[Completely Positive Cone: ${{\cC}^*_{n,d}}$]
	A tensor $\mathcal{X}\in {{\cS}_{n,d}}$ is completely positive if  $\exists ~ \bx_k\in \Re^{n}_{+}$ such that; $\mathcal{X}=\sum_{k=1}^{N} (\bx_{k})^{d}$.
	\begin{align}
	{{\cC}^*_{n,d}}&:=\left\{ \sum^{N}_{k=1}{(\bx_k) ^d}~\forall~ \bx_k\in {\Re^{n}_{+}}, N\in \mathbb{N} \right \}
	\end{align}
\end{definition}
\noindent 
The cones ${\cS}_{n,d}, {\cC}_{n,d} ~~and~~ {{\cC}_{n,d}}^*$ defined above, are all convex cones.
 we describe a special case for which the cones $\mathcal{S}_{n,d}^+$ and $\mathcal{C}_{n,d}$ coincide. The following theorem states the condition under which  a tensor having non-positive off-diagonal entries is copositivity.
\begin{theorem}
	Let $\mathcal{A} \in \mathcal{S}_{n,d}$ ($d$ is even) be an arbitrary tensor with all off-diagonal entries non-positive i.e. there exist $j$ and $k$ such that $i_{j} \ne i_{k}$ with $ \ a_{i_1i_2\cdots i_d} \le 0$ then $\mathcal{A} \in \mathcal{C}_{n,d}$ if and only if $\mathcal{A} \in \mathcal{S}_{n,d}^+$.
\end{theorem}
\noindent
\Pf
Let $\mathcal{A} \in \mathcal{C}_{n,d}$, then we have; 
\begin{align}\label{CopCondition}
f_\mathcal{A}(\bx)=\sum_{i_1,i_2 \cdots, i_d =1}^{n} a_{i_1i_2 \cdots i_d} x_{i_1} x_{i_2} \cdots x_{i_d}~\ge~0 ~ ~\forall~ \bx \in \Re^n_+
\end{align}
where all off-diagonal entries $a_{i_1i_2 \cdots i_d}  \le 0$. \\
To show that $\mathcal{A}$ is positive semidefinite tensor, we consider a polynomial form $f_{\mathcal{A}}(\bz)$ for all $\bz \in \Re^n$. Since $\bz$ has positive, negative and zero components. Therefore, to bifurcate among the non-negative and negative terms in the form $f_{\mathcal{A}}(\bz)$, we introduce few notations as follows; $\gamma_+(\bz):=\{k:z_{k} \ge 0\}$ and non-empty $ \gamma_-(\bz):=\{j:z_{j} < 0\}=\mathbb{N}_n\backslash\gamma_+(\bz)\ne \phi
$ where $\mathbb{N}_n=\{1,2,\cdots,n\}$, and
$|\sigma(i_j)|:=\text{number of times} ~i_j ~\text{occurs in the index of } a_{i_1i_2 \cdots i_d}$.
Using these notations we rewrite the polynomial form $f_{\mathcal{A}}(\bz)$ as follows:
\begin{small}
	\begin{align}\label{PsdCondition}
	f_{\mathcal{A}}(\bz)=\sum_{i_1 \cdots, i_d =1}^{n}  (\delta_{|\sigma(k)|\alpha_{k}})(\delta_{|\sigma(j)|\beta_{j}})a_{i_1i_2 \cdots i_d}\prod_{\begin{array}{c}
		k\in \gamma_+(\bz) \\ j\in \mathbb{N}_n\backslash\gamma_+(\bz)
		\end{array}}^{} {(\bz^{{\be}_{i_k}})}^{\alpha_{k}(i_k)}  {(\bz^{\be_{i_j}})}^{\beta_{j}(i_j)}
	\end{align}
\end{small}
\noindent
where, $\delta$ is Kronecker delta; and \\  
\begin{minipage}{0.3cm}
	\begin{align*}
	\alpha_k(i_k)=\begin{cases}
	\begin{matrix}
	& \alpha_{kk}  = \alpha_{k}  & if \ i_k=k \\
	&  0  & if \ i_k \ne k
	\end{matrix}
	\end{cases}
	\end{align*}
\end{minipage} 
and
\begin{minipage}{0.3cm}
		\begin{align*}
		 \beta_j(i_j)=\begin{cases}
		\begin{matrix}
		& \beta_{jj}  =\beta_{j} & if \ i_j=j \\
		& 0 & if \ i_j \ne j
		\end{matrix}
		\end{cases}
		\end{align*}
\end{minipage} \linebreak \\
Moreover,  $\sum_{k}^{}\alpha_{k}+\sum_{j}^{}\beta_{j}=d$ for all $\alpha_{k},\beta_{j} \in \mathbb{Z}_{d+1}$. From equation (\ref{PsdCondition}), it is clear that, all its terms corresponding to diagonal elements $a_{{(i)}^{d}}$ are given as follows:
\begin{align}\label{diagonal}
\varOmega^{(+)}=\underset{\ge 0}{\underbrace{\sum_{i=1}^{n} a_{(i)^d}(z_{i})^d}} ~~ \text{as both}~a_{(i)^d}~\text{and}~(z_{i})^d~\text{are positive}~\forall ~i.
\end{align}
We analyze those terms in (\ref{PsdCondition}) which correspond to off-diagonal elements of $\mathcal{A}$. The powers $\beta_j$ of negative components of $\bz$ are crucial in this analysis. Therefore we consider two cases as follows:
\begin{small}
	\begin{itemize}
		\item[(i)] If $\sum_{j}^{}\beta_j$ is even then all such terms are non-positive, that is;
		\begin{align}\label{negative}
	\varLambda^{(-)}=\underset{\le 0}{\underbrace{\bigg(\sum_{i_1 \cdots, i_d =1}^{n}  (\delta_{|\sigma(k)|\alpha_{k}})(\delta_{|\sigma(j)|\beta_{j}})a_{i_1\cdots i_d}\prod_{k\in \gamma_+(\bz)}^{} {(\bz^{\be_{i_k}})}^{\alpha_{k}(i_k)} \prod_{j\in \mathbb{N}_n\backslash\gamma_+(\bz)}^{} {(\bz^{\be_{i_j}})}^{2\beta'_{j}(i_j)}\bigg)} }
		\end{align} 
		as $a_{i_1i_2 \cdots i_d}\le0$ and $\prod_{k\in \gamma_+(\bz)}^{} {(\bz^{e_{i_k}})}^{\alpha_{k}(i_k)} \prod_{j\in \mathbb{N}_n\backslash\gamma_+(\bz)}^{} {(\bz^{e_{i_j}})}^{2\beta'_{j}(i_j)} \ge 0$ therefore sum of non-positive terms is non-positive again.
		\item[(ii)] If $\sum_{j}^{}\beta_j$ is odd then all such terms are non-negative, that is;
		\begin{align}\label{positive}
	\varGamma^{(+)}=\underset{\ge 0}{\underbrace{\bigg(\sum_{i_1 \cdots, i_d =1}^{n}  (\delta_{|\sigma(k)|\alpha_{k}})(\delta_{|\sigma(j)|\beta_{j}})a_{i_1\cdots i_d}\prod_{k\in \gamma_+(\bz)}^{} {(\bz^{e_{i_k}})}^{\alpha_{k}(i_k)} \prod_{j\in \mathbb{N}_n\backslash\gamma_+(\bz)}^{} {(\bz^{e_{i_j}})}^{\beta_{j}(i_j)}\bigg)}}	
		\end{align} 
		as two factors $a_{i_1i_2 \cdots i_d}$ and $\prod_{j\in \mathbb{N}_n\backslash\gamma_+(\bz)}^{} {(\bz^{\be_{i_j}})}^{\beta_{j}(i_j)}$ are non-positive, and the third one $\prod_{k\in \gamma_+(\bz)}^{} {(\bz^{e_{i_k}})}^{\alpha_{k}(i_k)}$ is non-negative therefore the product of these three factors is non-negative. Hence the sum of non-negative terms is non-negative again.
	\end{itemize} 
\end{small}
Summarizing equation (\ref{PsdCondition}) by incorporating the information given in equations (\ref{diagonal}),(\ref{positive})and (\ref{negative}) we have: 
\begin{align}\label{=13}
f_\mathcal{A}(\bz)&=\varOmega^{(+)}+\varLambda^{(-)}+\varGamma^{(+)}\\
&=\varOmega-\varLambda+\varGamma ~~\text{where}~\varLambda^{(-)}=-\varLambda~~\text{and}~\varOmega,\varLambda, \varGamma \in \Re_+
\end{align}
We define a mapping $\phi:\Re^n \rightarrow\Re^n_+$ as follows;
\begin{align*}
\phi(\bz)=\begin{pmatrix}
a_{11} & 0 & \cdots & 0 \\ 0 & a_{22} & \cdots & 0\\ \vdots & \vdots & \ddots & \vdots \\ 0 & 0 & \cdots & a_{nn}
\end{pmatrix} \begin{pmatrix}
z_1 \\ z_2 \\ \vdots \\ z_n
\end{pmatrix} 
~~where~~a_{ii}=\bigg \lbrace \begin{array}{cc}
~ ~1 & if ~i \in \gamma_+(\bz) \\
-1 & if ~i \in \gamma_-(\bz)
\end{array}
\end{align*}
thus corresponding to each $\bz\in\Re^n$ there exist unique $\phi(\bz) \in \Re^n_+$, 
therefore,  rewriting (\ref{CopCondition}) as a sum of non-negative and non-positive forms as follows: 
\begin{align}\label{19}
f_\mathcal{A}(\phi(\bz))=\sum_{i=1}^{n} \underset{\ge 0}{\underbrace{{a_{(i)^d}}}} ~ \underset{\ge 0}{\underbrace{{(\phi(z_i))^d}}} + \sum_{i_1,i_2 \cdots, i_d =1}^{n} \underset{\le 0}{\underbrace{a_{i_1i_2 \cdots i_d}}} ~  \underset{\ge 0}{\underbrace{\phi(z_{i_1}) \phi(z_{i_2}) \cdots \phi(z_{i_d})}}
\end{align} 
for even order $d$ the first and second terms in (\ref{19}) are non-negative and  non-positive respectively, therefore we have;
\begin{align}\label{14}
\varOmega^{(+)}=\sum_{i+1}^{n} & a_{(i)^d}{\phi(z_i)}^d  ~ ~ ~  and ~ ~ ~~\varPsi^{(-)}= \sum_{i_1,i_2 \cdots, i_d =1}^{n} a_{i_1i_2 \cdots i_d} \phi(z_{i_1})\phi(z_{i_2}) \cdots \phi(z_{i_d})
\end{align}
Clearly all the terms in $\varPsi^{(-)}$ are non-positive, that is, $\varPsi^{(-)}=-\varPsi$ where $\varPsi \in \Re_+$, thus we have; 
\begin{align*}
& f_\mathcal{A}(x)=\varOmega^{(+)}+\varPsi^{(-)}\\
\implies & f_\mathcal{A}(x)=\varOmega-\varPsi ~~\text{where}~~\varOmega,\varPsi \in \Re_+ \\
\implies & f_\mathcal{A}(x)=\varOmega-(\varLambda+\varGamma) ~~\because~\varPsi=\varLambda+\varGamma
\end{align*} 
Clearly, we have;
\begin{align*}
\varLambda-\varGamma \le & ~ ~\varLambda+\varGamma \\
\implies -(\varLambda - \varGamma )\ge & -(\varLambda+\varGamma) \\
\implies \varOmega-(\varLambda-\varGamma) \ge & ~\varOmega-(\varLambda+\varGamma) ~~\because~\varOmega\in \Re_+\\
\implies f_\mathcal{A}(\bz)=\varOmega-(\varLambda-\varGamma )\ge & \varOmega-(\varLambda+\varGamma) =f_\mathcal{A}(\bx) \ge  0 ~~ \text{for}~\bx\in \Re^n_+~\text{and}~\bz\in \Re^n
\end{align*}
Thus,$f_\mathcal{A}(\bz) \ge 0~~\forall~\bz\in \Re^n$ which means, $\mathcal{A}\in \mathcal{S}_{n,d}^+$. \\
Converse is straight forward as all positive semidefinite tensors are copositive as well. 
\eop \\
Here, we would like to remark that the above theorem has a significant importance for the case $d=2$, since it gives a special class consisting of polynomial time solvable problems as a sub-class of  an in general NP-hard problem. However, for the case when $d\ge4$ this property is not true due to the fact that determining if a tensor is PSD is an NP-hard problem \cite{Hillar:2013:MTP:2555516.2512329}. 
\section{Approximation Hierarchies for Copositive Cone of Tensors}
The copositive cone of tensors is not computationally tractable. To solve optimization programs which involves such cones, we replace the copositive (completely positive) cone by its approximation hierarchies which yields an approximate optimal value of the original program. In this section we present, several inner and outer approximation for the copositive cone of tensors.
\subsection{Inner Approximation Hierarchies for Copositive Cone of Tensors}
In the first part of this section, we present polynomial based inner approximation hierarchies  $\mathcal{C}_{n,d}^{(r)}$ and $\mathcal{K}_{n,d}^{(r)}$ for copositive cone of tensors. Secondly, we present the inner approximations $\mathcal{I}_{n,d}^{\mathcal{P}}$ for copositive cone based on simplicial partition $\mathcal{P}$ of the simplex $\Delta$. In the last part,  the containment relations among the above mentioned inner approximations are discussed. 
\subsubsection{Inner Approximations Based on Polynomial Conditions}
Since for any polynomial $f_\mathcal{A}(\bx)$ of degree $d$ in $n$ variables as defined in (\ref{2}), by virtue of the fact that for any $\bx\in\Re^n_+$ we may write ${\bx}=\by \circ \by$ for some $\by\in\Re^n$, where $\circ$ indicates the component wise product; thus the copositivity conditions for a tensor $\mathcal{A}$ can be represented as follows;
\begin{align}
f_\mathcal{A}{(\by)} &=\bigg\langle \mathcal{A}, {\mathcal{T}_d(\by \circ \by)} \bigg\rangle \\
 &=  \bigg(\sum_{i_1, \ldots, i_d =1}^n a_{{i_1}{i_2} \ldots i_d}{{ y^2_{i_1}}} { {y^2_{i_2}}}\cdots { {y^2_{i_d}}}\bigg) \ge 0 ~~ \forall~ \by\in\Re^n
\end{align} 
In the subsequent part of this section, we adapt the idea from Parrilo \cite{Parrilo:2013:CAG:2465506.2466575},  De Klerk and Pasechnik \cite{deklerk2002} for higher degree sequence of polynomials  $\lbrace P^{(r)}(\by)\rbrace _{r\in \mathbb{N}_{0}}$ which is stated as under : 
\begin{align}\label{sos}
P^{(r)}\big(\by\big)=f_\mathcal{A}{(\by)} \bigg(\sum_{k=1}^{n}y^2_k\bigg)^r~~ for~all ~~\by\in\Re^n
\end{align}
\noindent
Based on the above representation of polynomials, and by virtue of algebraic geometry Pablo A. Parrilo \cite{Parrilo:2013:CAG:2465506.2466575} defined a sufficient condition for the polynomial to be non-negative everywhere, which further culminates an idea to define an approximation for the copositive cone  known as the Parrilo cone. The tensor analogue of Parrilo cone denoted by $\mathcal{K}_{n,d}^{(r)}$ which consists of those tensors $\mathcal{A} \in {\mathcal{S}}_{n,d}$ 
for which the associated polynomial in (\ref{sos}) has sums of squares decomposition. That is;
\[\mathcal{K}_{n,d}^{(r)}=\bigg\lbrace \mathcal{A} \in {\mathcal{S}}_{n,d}: P^{(r)}(\by)=f_\mathcal{A}{(\by)} \left(\sum_{k=1}^{n}y^2_k\right)^{(r)} \text{ has an SOS decomposition} \bigg\rbrace \] 
Clearly, $\mathcal{K}_{n,d}^{(r)} \subseteq \mathcal{K}_{n,d}^{(r+1)}$ for $r\in\mathbb{N}_{0}$, since:
\begin{align*}
P^{(r+1)}(\by)&=f_\mathcal{A}{(\by)} \left(\sum_{k=1}^{n}y^2_k\right)^{(r+1)}\\ &=\left(\sum_{k=1}^{n}y^2_k\right) \left(f_\mathcal{A}{(\by)}\left(\sum_{k=1}^{n}y^2_k\right)^{(r)}\right)\\
&= \left(\sum_{k=1}^{n}y^2_k\right)P^{(r)}(\by).
\end{align*}
For $r=0$ we have,
\begin{align*}
\mathcal{K}_{n,d}^{(0)}&=\bigg\lbrace \mathcal{A} \in {\mathcal{S}}_{n,d}: P^{(0)}(\by)=f_\mathcal{A}{(\by)} \left(\sum_{k=1}^{n}y^2_k\right)^{(0)} \text{ allows an SOS decomposition} \bigg\rbrace  \\
&=\bigg\lbrace \mathcal{A} \in {\mathcal{S}}_{n,d}: f_\mathcal{A}{(\by)} \text{ allows an SOS decomposition} \bigg\rbrace
\end{align*}  
As we know that, if $f_\mathcal{A}{(\by)}$ allows an SOS decomposition then its associated tensor $\mathcal{A} \in \mathcal{S}_{n,d}^+$, however converse is not true in general (see counter Example 4.5, \cite{Parrilo00structuredsemidefinite}). Clearly the following inclusion relation holds;
\[\mathcal{K}_{n,d}^0\subset \mathcal{K}_{n,d}^1\subset \cdots \subset \mathcal{K}_{n,d}^{\infty}=\mathcal{C}_{n,d} \]
where $\mathcal{S}_{n,d}$ denotes the cone of positive semidefinite tensors.\\
Another, sufficient condition for a polynomial to be non-negative, defined by Bomze and De Klerk \cite{bomze2002solving}, cf. also De Klerk and Pasechnik \cite{deklerk2002} which exploits the coefficients. The tensor analogue of this condition leads us to define  an approximation cone $\mathcal{C}_{n,d}^{(r)}$ for the copositive cone which consists of those tensors $\mathcal{A}\in \mathcal{S}_{n,d}$ for which the associated polynomial in (\ref{sos}) has non-negative coefficients. That is; 
\[\mathcal{C}_{n,d}^{(r)}=\bigg\lbrace \mathcal{A} \in {\mathcal{S}}_{n,d}: P^{(r)}(\by)=f_\mathcal{A}{(\by)} \left(\sum_{k=1}^{n}y^2_k\right)^{(r)} \text{ has non-negative coefficients} \bigg\rbrace\] 
 For an arbitrary $\mathcal{A} \in \mathcal{C}_{n,d}^{(r)}$ and due to the fact that $P^{(r+1)}(\by)= \big(\sum_{k=1}^{n}y^2_k\big)P^{(r)}(\by)$,  it is straightforward to see that $\mathcal{C}_{n,d}^{(r)} \subseteq \mathcal{C}_{n,d}^{(r+1)}$. Moreover, if $\mathcal{A} \in \mathcal{C}_{n,d}^{(r)}$ then definitely its associated polynomial $P^{(r)}(\by)$ has non-negative coefficients which further allows to have an SOS decomposition, thus $\mathcal{A} \in  \mathcal{K}_{n,d}^{(r)}$ for each $r\in\mathbb{N}_{0}$.  Therefore we have; 
\[\mathcal{C}_{n,d}^{(r)} \subseteq \mathcal{K}_{n,d}^{(r)} \ \forall \ r\in\mathbb{N}_{0} \]
Both cones $\mathcal{C}_{n,d}^{(r)}$ and $\mathcal{K}_{n,d}^{(r)}$ are inner approximations of $\mathcal{C}_{n,d}$, that is 
\[int\big(\mathcal{C}_{n,d}\big)\subseteq \bigcup_{r\in \mathbb{N}_{0}} \mathcal{K}_{n,d}^{(r)} \subseteq \mathcal{C}_{n,d} ~ ~ \text{and} ~ ~ int\big(\mathcal{C}_{n,d}\big)\subseteq \bigcup_{r\in \mathbb{N}_{0}} \mathcal{C}_{n,d}^{(r)} \subseteq \mathcal{C}_{n,d}\]
We discuss the procedure to calculate the coefficients of polynomials $P^{(r)}(\by)$ of degree $2(d+r)$. 
Generalizing the characterizations  given by Bomze and De Klerk in \cite{bomze2002solving} for the tensors of order $d$ and dimension $n$. For an arbitrary $\balpha \in \Re^n$ we define the multinomial coefficients as follows: 
\begin{align}\label{c(m)}
c(\balpha)=\begin{cases}
\begin{matrix}
\frac{\|\balpha\|_{1}!}{\prod_{i}^{n}(\alpha_i!)}  
& if~\balpha\in\mathbb{N}^n_0 \\
0  & \ \ if~\balpha\in \Re^n \backslash \mathbb{N}^n_0
\end{matrix}
\end{cases}
\end{align}
where for any integer $k\in\mathbb{Z}$, $k!$ denotes the factorial of $k$. Expanding the polynomial $P^{(r)}(y)$ in (\ref{sos}) by using the multinomial law we have the following:
\begin{align}\label{22}
P^{(r)}({\by})&=f_\mathcal{A}{({\by})} \left(\sum_{k=1}^{n}y^2_k\right)^r\\
&=\bigg(\sum_{i_1, \ldots, i_d =1}^{n} a_{i_1i_2 \ldots i_d}{{ y^2_{i_1}}} { {y^2_{i_2}}}\cdots { {y^2_{i_d}}}\bigg)\bigg(\sum_{\balpha\in\mathbb{I}^n(r)}^{}c(\balpha) \by^{2\alpha}\bigg)\\
&=\bigg(\sum_{i_1, \ldots, i_d =1}^n a_{i_1 \ldots i_d}{{ \by}^{2\be_{i_1}}} { {\by}^{2\be_{i_2}}}\cdots { {\by}^{2\be_{i_d}}}\bigg)\bigg(\sum_{\balpha\in\mathbb{I}^n(r)}^{}c(\balpha) \by^{2\balpha}\bigg)\\
&=\sum_{\balpha\in\mathbb{I}^n(r)}^{} \sum_{i_1, \ldots, i_d =1}^n c(\balpha) a_{i_1 \ldots i_d} \by^{2(\balpha+\be_{i_1}+\be_{i_2}+\cdots+\be_{i_d})}\\
&=\sum_{\balpha\in\mathbb{I}^n(r)}^{} \bigg(\sum_{i_1, \ldots, i_d =1}^{n} c(\balpha) a_{i_1 \ldots i_d}\bigg) \by^{2(\balpha+\be_{i_1}+\be_{i_2}+\cdots+\be_{i_d})}
\end{align} 
let us assume that, $\btheta=\balpha+\be_{i_1}+\be_{i_2}+\cdots+\be_{i_d}$, and abbreviating \[\btheta\big(i_1,i_2,\cdots,i_d\big)=\btheta-(\be_{i_1}+\be_{i_2}+\cdots+\be_{i_d})\]
taking $s=r+d$, from the last identity in (\ref{22}) it follows that;
\begin{align}\label{multicoeff}
P^{(r)}(\by)=\sum_{\btheta\in\mathbb{I}^n(s)}^{} \bigg(\sum_{i_1, \ldots, i_d =1}^{n} c\big(\btheta(i_1,i_2,\cdots,i_d)\big) a_{i_1i_2 \ldots i_d}\bigg) \by^{2(\btheta)}
\end{align}
denoting the coefficients in (\ref{multicoeff}) by $\mathcal{A}_{\btheta}$, that is we have;
\begin{align}\label{Am}
\mathcal{A}_{\btheta}=\sum_{i_1, \ldots, i_d =1}^n c\big(\btheta(i_1,i_2,\cdots,i_d)\big) a_{i_1i_2 \ldots i_d}
\end{align} 
The procedure for finding the coefficients $c\big(\btheta(i_1,i_2,\cdots,i_d)\big)$ of multinomial in (\ref{Am}) is described in the following proposition. 
\begin{proposition}
	Let $c\big(\btheta(i_1,i_2,\cdots,i_d)\big)$ be the coefficient of multinomial $P^{(r)}(\by)$ in (\ref{multicoeff}), then we have;
	\begin{align*}
	c\big(\btheta(i_1,\cdots,i_d)\big)=\bigg\lbrace  
	\begin{matrix}
	c\big(\btheta-d\be_{i}\big);& i=i_j~:~\forall~j\in \mathbb{N}_{d}~and~i\in\mathbb{N}_{n}\\
	0; & if ~~ i=i_1=\cdots=i_k \ne i_{k+1}\ne \cdots \ne i_{d}\\
	c\big(\btheta-(\be_{i_1}+\cdots+\be_{i_d})\big);&if ~~ i_1 \ne i_2 \ne  \cdots \ne i_{d}\\	
	\end{matrix}
	\end{align*}
\end{proposition}
\noindent
\Pf
Since, $\btheta=\alpha+\be_{i_1}+\be_{i_2}+\cdots+\be_{i_d}$ therefore $\|\btheta\|_{1}=r+d$. By (\ref{c(m)}) we see that, when index $i$ is repeated $d-times$ then we have, $c\big(\btheta(i,i,\cdots,i)\big)\ne0~~as~~\omega_i > d$. However, when some index $i$ is repeated $k-times$ where $1<k<d$, then $\omega_i=r_{i}+k$, so there exist an active index $j\ne i$ such that  $m_j\le0$, which implies that, such coefficients are always zero; that is 
\begin{align*}
c\big(\btheta(i_1,\cdots,i_d)\big)=0; ~~ if ~~ i=i_1=\cdots=i_k \ne i_{k+1}\ne \cdots \ne i_{d}
\end{align*}
 we consider the last case when, $i_1 \ne i_2 \ne  \cdots \ne i_{d}$. In such case we have, $m_{j}=r_{j}-1~~\forall j \in \mathbb{N}_d$, thus we have;
\begin{align*}
c\big(\btheta(i_1,\cdots,i_d)\big)= c\big(\btheta-(\be_{i_1}+\cdots+\be_{i_d})\big).
\end{align*}
\eop
Analogous to matrix case, a vector extracted from all diagonal entries of $\mathcal{A}$ is denote by $diag(\mathcal{A})\in\Re^n$, i.e. \[diag(\mathcal{A})=\begin{bmatrix}
a_{(i_1)^d}\\a_{(i_2)^d}\\ \vdots \\a_{(i_n)^d}
\end{bmatrix} \] 
However, $Diag(\btheta)$ is a diagonal tensor of order $d$ and dimension $n$ having components of $\btheta$ at its diagonals.
By using the definition (\ref{c(m)}) of $c(\btheta)$ the representation of $\mathcal{A}_{\btheta}$ given in  (\ref{multicoeff}) is simplified considerably in the following theorem.
\begin{theorem}\label{propo1}
	Let $\mathcal{A}\in \mathcal{S}_{n,d}$ be any $d^{th}$-order $n$-dimensional symmetric tensor then, for $ \btheta\in \mathbb{I}^n(s)$ we have;
	\begin{small}
		\begin{align*}\label{propo}
		\mathcal{A}_{\btheta}=&\frac{c\big(\btheta\big)}{s(s-1)(s-2)\cdots(s-(d-1))} \bigg[\bigg\langle \mathcal{A},\btheta^d \bigg\rangle + \\
		&\sum_{k=1}^{(d-1)}(-1)^{k} \bigg(\sum_{ \prod_{j=1}^{k}{\theta_{t_j}}\in \{\sigma (i_1i_2\cdots i_k) : \{i_j\}_{j=1}^{k}\subseteq \mathbb{N}_{d-1} \}}^{}  \big(\prod_{j=1}^{k}{\theta_{t_j}}\big)\bigg)\bigg \langle \mathcal{A},Diag~\underset{(d-k)-times}{\underbrace{ \big(\btheta\circ\cdots\circ\btheta\big)}}\bigg \rangle \bigg] 
		\end{align*}
	\end{small}
\end{theorem}
\noindent
\Pf  
For $i_1=i_2=\cdots=i_d=i$ the coefficients $c\big(\btheta(i,i,\cdots,i)\big)$ are zero if $\omega_i<d$, however for $i_1 \ne i_2\ne\cdots\ne i_d$ the coefficients $c\big(\btheta(i_1,i_2,\cdots,i_d)\big)=0$ if $\prod \omega_i=0$. Thus the nonzero coefficients in (\ref{multicoeff}) occurs only for some $(i_1,i_2,\cdots,i_d)$ tuples depending upon $\btheta$. Therefore after some straightforward calculations, we have the following simplification of  (\ref{Am});
\begin{small}
	\begin{align*}
	\mathcal{A}_{\btheta}&=\sum_{i_1, \ldots, i_d =1}^{n} c\big(\btheta(i_1,i_2,\cdots,i_d)\big) a_{i_1i_2 \ldots i_d}\\
	&=\sum_{{i_1,i_2, \ldots, i_d} =1 }^{n}{\frac{\|\btheta-(\theta_1 \be_{i_1}+\cdots+\theta_d \be_{i_d})\|_{1}!} {\prod_{i}^{n}(\omega_i-{\theta_i})!}} a_{i_1i_2\ldots i_d};~~with~ \sum_{i}^{}\theta_i=d ~\forall~\theta_i \in \mathbb{Z}_{d+1}\\
	&=\sum_{i=1}^{n}{\frac{\|\btheta-{d}\be_{i}\|_{1}!} {(\omega_i-{d})!\prod_{k \in \mathbb{N}_{n}\backslash \{i\}}^{}(\omega_k)!}} a_{(i)^d} ~+\\
	& \sum_{ \begin{tiny}
		\begin{array}{c}
			i_1,i_{j}=1 \\
		i_j \in \lbrace i_j:i_j=i_k~\forall~ j,k \in \mathbb{N}_{d} \backslash \{1\} \rbrace 
		\end{array}
		\end{tiny}}^{n}\frac{\|\btheta-\be_{i_1}-(d-1)\be_{i_j}\|_{1}!} {(\omega_{i_1}-{1})!{(\omega_{i_j}-{(d-1)})!\prod_{k \in \mathbb{N}_{n}\backslash \{1,j\}}^{}(\omega_k)!}} a_{i_1(i_j)^{(d-1)} }  + \\
	&\sum_{ \begin{tiny}
		\begin{array}{c}
		i_1,i_2,i_j=1 \\
		i_j \in \lbrace i_j : i_j=i_k~\forall~ j,k \in \mathbb{N}_{d} \backslash \{1,2\} \rbrace 
		\end{array}
		\end{tiny} }^{n}\frac{\|\btheta-\be_{i_1}-\be_{i_2}-(d-2)\be_{i_j}\|_{1}!} {(\omega_{i_1}-{1})!(\omega_{i_2}-{1})!{(\omega_{i_j}-{(d-2)})!\prod_{k \in \mathbb{N}_{n}\backslash \{1,2,j\}}^{}(\omega_k)!}} a_{i_1i_{2}(i_j)^{(d-2)} }+\cdots  \\
	&+\sum_{ \begin{tiny}
		\begin{array}{c}
		i_1,i_2,\cdots,i_{(d-1)},i_j=1 \\
		i_j \in \lbrace i_j : i_j=i_k~\forall~ j,k \in \mathbb{N}_{d} \backslash \{1,2,\cdots,d-1\} \rbrace 
		\end{array}
		\end{tiny}}^{n}\frac{\|\btheta-\be_{i_1}-\cdots-\be_{i_d}\|_{1}!} {\prod_{k=1}^{d} (\omega_{k}-{1})!\prod_{k \in \mathbb{N}_{n}\backslash \{1,2,\cdots,d-1\} }^{}(\omega_k)!} a_{i_1i_2 \ldots i_d} 
	\end{align*}
	using $\|\btheta\|_{1}=s$ and by virtue of the fact that the difference between positive numbers, it implies that, $\|\btheta-(\theta_1 \be_{i_1}+\cdots+\theta_d \be_{i_d})\|_{1}=(s-d)$, which further simplifies the previous equation as follows:
	\begin{align*}
	\mathcal{A}_{\btheta}&=\frac{c(\btheta)}{s(s-1)(s-2)\cdots(s-(d-1))} \Bigg[\sum_{i=1}^{n}{\bigg(\prod_{k=0}^{d-1} (\omega_i-k)}\bigg) a_{(i)^d} + \\
	&~\ \ \ \ \ \ \ \ \ \ \ \ \ \ \ \ \ \ ~ \sum_{ \begin{tiny}
		\begin{array}{c}
		i_1,i_{j}=1 \\
		i_j \in \lbrace i_j:i_j=i_k~\forall~ j,k \in \mathbb{N}_{d} \backslash \{1\} \rbrace 
		\end{array}
		\end{tiny}}^{n}{\bigg(\omega_{i_1} \prod_{k=0}^{d-2} (\omega_{i_j}-k)\bigg)} a_{i_{1}(i_j)^{(d-1)}}+\\
	&~~ \sum_{ \begin{tiny}
		\begin{array}{c}
		i_1,i_2,i_j=1 \\
		i_j \in \lbrace i_j : i_j=i_k~\forall~ j,k \in \mathbb{N}_{d} \backslash \{1,2\} \rbrace 
		\end{array}
		\end{tiny} }^{n}{ \bigg( \omega_{i_1} \omega_{i_2}\prod_{k=0}^{d-3} (\omega_{i_j}-k)\bigg) {a_{{i_1}i_{2}(i_j)^{(d-2)}}} } + \cdots +\\ 
	&~\ \ \ \ \ \ \ \ \ \ \ \ \ \ \ \ \ \  \ \  \sum_{ \begin{tiny}
		\begin{array}{c}
		i_1,i_2,\cdots,i_{(d-1)},i_j=1 \\
		i_j \in \lbrace i_j : i_j=i_k~\forall~ j,k \in \mathbb{N}_{d} \backslash \{1,2,\cdots,d-1\} \rbrace 
		\end{array}
		\end{tiny}}^{n} \bigg(\prod_{k=1}^{d} \omega_{i_k} \bigg) {a_{i_1i_2 \ldots i_d} } \Bigg]\\ 
	\end{align*}
	since by multi-coefficient formula (\ref*{c(m)}) it is evident that $c\big(\btheta(i_1,i_2,\cdots,i_d)\big)=0$, as some of the  components of $\btheta(i_1,i_2,\cdots,i_d)$ can be negative, or the product of its components is zero. Therefore, the coefficients of  $a_{i_{1},(i_j)^{(d-1)}},a_{i_{1},i_2,(i_j)^{(d-2)}},\cdots, a_{i_{1},i_2,\cdots,(i_{(d-1)})} $ vanish. Thus, we have the following simplified equation: 
	\begin{align*}
	\mathcal{A}_{\btheta}&=\frac{c(\btheta)}{s(s-1)(s-2)\cdots(s-(d-1))} \Bigg[\sum_{i=1}^{n}{\bigg(\prod_{k=0}^{d-1} (\omega_i-k)}\bigg) a_{(i)^d} + \\
	& \ \ \ \ \ \ \ \ \ \ \ \ \ \ \ \ \ \ \ \ \ \ \sum_{ \begin{tiny}
		\begin{array}{c}
		i_1,i_2,\cdots,i_{(d-1)},i_j=1 \\
		i_j \in \lbrace i_j : i_j=i_k~\forall~ j,k \in \mathbb{N}_{d} \backslash \{1,2,\cdots,d-1\} \rbrace 
		\end{array}
		\end{tiny}}^{n} \bigg(\prod_{k=1}^{d} \omega_{i_k} \bigg) {a_{i_1i_2 \ldots i_d} } \Bigg]\\ 
	&=\frac{c\big(\btheta\big)}{s(s-1)(s-2)\cdots(s-(d-1))}\bigg[\sum_{i}^{}\bigg( \omega_i\big(\omega_i-1\big)\cdots\big(\omega_i-(d-1)\big)\bigg) a_{(i)^d} +\\
	& \ \ \ \ \ \ \ \ \ \ \ \ \ \ \ \ \ \  \sum_{ \begin{tiny}
		\begin{array}{c}
	i_1,i_2,\cdots,i_{(d-1)},i_j=1 \\
	i_j \in \lbrace i_j : i_j=i_k~\forall~ j,k \in \mathbb{N}_{d} \backslash \{1,2,\cdots,d-1\} \rbrace 
		\end{array}
		\end{tiny}}^{n}\bigg(\omega_{i_1}\omega_{i_2}\cdots \omega_{i_d}\bigg)a_{i_1i_2 \ldots i_d}\bigg]\\
	&=\frac{c\big(\btheta\big)}{s(s-1)(s-2)\cdots(s-(d-1))}\bigg[\sum_{
		\begin{tiny}
		\begin{array}{c}
		i_1,i_2,\cdots,i_{(d-1)},i_j=1 \\
		i_j \in \lbrace i_j : i_j=i_k~\forall~ j,k \in \mathbb{N}_{d} \backslash \{1,2,\cdots,d-1\} \rbrace 
		\end{array}
		\end{tiny}}^{n}\big(\omega_{i_1}\cdots \omega_{i_d}\big)a_{i_1 \cdots i_d}+\\
	&\sum_{i=1}^{n}\big(\omega_i^d\big)a_{(i)^d} + \sum_{i=1}^{n} \bigg(\sum_{k=1}^{(d-1)}(-1)^{k} \bigg(\sum_{ \big(\prod_{j=1}^{k}{\theta_{t_j}}\big)\in \{\sigma (i_1i_2\cdots i_k) : \{i_j\}_{j=1}^{k}\subseteq \mathbb{N}_{d-1} \}}^{}  \big(\prod_{j=1}^{k}{\theta_{t_j}}\big)\bigg){\omega_i}^{(d-k)} \bigg){a_{(i)^d}}\bigg]\\
	&=\frac{c\big(\btheta\big)}{s(s-1)(s-2)\cdots(s-(d-1))}\bigg[\bigg\langle \mathcal{A},\btheta^d \bigg\rangle +\\
	&~~ \bigg(\sum_{k=1}^{(d-1)}(-1)^{k} \bigg(\sum_{ {(\prod_{j=1}^{k}{\theta_{t_j}})}\in \{\sigma (i_1i_2\cdots i_k) : \{i_j\}_{j=1}^{k}\subseteq \mathbb{N}_{d-1} \}}^{}  \big(\prod_{j=1}^{k}{\theta_{t_j}}\big)\bigg) \bigg\langle \mathcal{A},Diag\underset{(d-k)-times}{\underbrace{ \big(\btheta\circ\btheta\circ\cdots\circ\btheta\big)}}\bigg\rangle\bigg)\bigg]
	\end{align*}
\end{small}
\eop
For the sake of brevity, we denote the coefficients of inner product $\bigg\langle \mathcal{A},Diag\underset{(d-k)-times}{\underbrace{ \big(\btheta\circ\btheta\circ\cdots\btheta\big)}} \bigg\rangle$ in the last identity by $\beta_k$, that is; 
\[\beta_k=\bigg(\sum_{ {(\prod_{j=1}^{k}{\theta_{t_j}})}\in \lbrace\sigma (i_1i_2\cdots i_k): \{i_j\}_{j=1}^{k}\subseteq \mathbb{N}_{d-1} \rbrace }^{}  \big(\prod_{j=1}^{k}{\theta_{t_j}}\big)\bigg)~~where~~k=1,2,\cdots,(d-1) \] 
Notice from the above that $P^{(r)}(\by)$ have non-negative coefficient whenever $\mathcal{A}_{\btheta} \ge 0$. Thus we arrive at the following theorem. 
\begin{theorem}
	For an arbitrary $\btheta\in\Re^n$ and $r\in \mathbb{N}_{0}$, the cone $\mathcal{C}_{n,d}^{(r)}$ is defined as:
	\begin{small}
		\begin{align*}
	\mathcal{C}_{n,d}^{(r)}=\left\lbrace \mathcal{A} \in {\mathcal{S}}_{n,d}:\bigg\langle \mathcal{A},\bigg(\btheta^d+\sum_{k=1}^{(d-1)}(-1)^{k} \beta_k Diag~\big(\underset{(d-k)-times}{\underbrace{ \btheta\circ\cdots\circ \btheta}} \big)\bigg)\bigg\rangle\ge 0~~\forall~ \btheta\in \mathbb{I}^n(s)\right\rbrace
	\end{align*}
	\end{small}
\end{theorem}
\noindent
\Pf It is an immediate consequence of  (\ref{multicoeff}) and Proposition \ref{propo1}. 
\eop \\
 for $r=0$ we have;
\begin{small}
	\begin{align*}
	\mathcal{C}_{n,d}^{(0)}=\left\lbrace \mathcal{A} \in {\mathcal{S}}_{n,d}:\bigg\langle \mathcal{A} \ ,\bigg(\btheta^d+\sum_{k=1}^{(d-1)}(-1)^{k} \beta_k Diag~\big(\underset{(d-k)-times}{\underbrace{ \btheta\circ\cdots\circ \btheta}} \big)\bigg)\bigg\rangle\ge 0~~\forall~ \btheta\in \mathbb{I}^n(d)\right\rbrace
	\end{align*}
\end{small}
Clearly, if $\mathcal{A} \in \mathcal{C}_{n,d}^{(0)}$ then the entries of $\mathcal{A}$ must be non-negative. Thus $\mathcal{A} \in \mathcal{N}_{n,d}$; that is, 
\[\mathcal{N}_{n,d}=\mathcal{C}_{n,d}^0\subset\mathcal{C}_{n,d}^1\subset \cdots \subset \mathcal{C}_{n,d}^{\infty} =\mathcal{C}_{n,d} \] 
\subsubsection{Inner Approximations Based on Simplicial Partition}
Let $\|\cdot\|_{1}$ denote 1-norm on $\Re^n_+.$ The set $\Delta=\{\bx\in\Re^n_+:\|x\|_1=1\}$ is known as the standard
simplex. Any tensor $\mathcal{A}$ is copositive if and only if $\langle \mathcal{A}, \mathcal{T}_{d}(\bx) \rangle \ge 0 ~~\forall~~\bx \in \Delta$. In the subsequent part of this section we state the conditions for non-negativity of the polynomial $f_\mathcal{A}(\bx)$ over a simplex. An appropriate way to express polynomials over a simplex $\Delta=conv\{\bv_1,\bv_2,\cdots,\bv_n\}$ is by using barycentric coordinates; that is, 
$\bx=\sum_{i=1}^{n}{\lambda_i \bv_i}$ where $\sum_{i=1}^{n}{\lambda_i }=1$
for $\lambda_1,\lambda_2,\cdots,\lambda_n\in\Re $ and $\bx \in \Delta$. The representation of polynomial form in barycentric coordinates is given as follows:
\begin{align*}
f_\mathcal{A}{({\bx})} &=\bigg\langle \mathcal{A} \ , {\mathcal{T}_{d}\bigg(\sum_{i=1}^{n}\lambda_i \bv_i\bigg)} \bigg\rangle \\
&=  \sum_{i_1, \ldots, i_d =1}^n {\lambda_{i_1}\lambda_{i_2} \cdots \lambda_{i_d}} \bigg\langle \mathcal{A} \ , \bv_{i_1}\otimes \bv_{i_2}\otimes \cdots\otimes \bv_{i_d} \bigg\rangle
\end{align*}
For the non-negative $\lambda$ and basis vectors of the simplex $\Delta$, we state the following lemma.
\begin{lemma}\label{lemma1}
	Let $\Delta=conv\{\bv_1,\bv_2,\cdots,\bv_n\}$ be a simplex. If $\big\langle \mathcal{A}, \bv_{i_1}\otimes \cdots\otimes \bv_{i_d}\big\rangle \ge 0$ for all $i_1, \cdots, i_d \in \{1,2,\cdots,n\}$ then $\bigg\langle \mathcal{A},\mathcal{T}_d(\bx) \bigg\rangle \ge 0 ~~\forall~\bx\in\Delta.$
\end{lemma}
\noindent
\Pf
Let $\Delta=conv\{\bv_1,\bv_2,\cdots,\bv_n\}$ then each $x\in\Delta$, can be expressed as;
\[\bx=\sum_{i=1}^{n}{\lambda_i \bv_i} ~~with~~\sum_{i=1}^{n}{\lambda_i }=1\]
 taking, 
\begin{align*}
\bigg\langle \mathcal{A},\mathcal{T}_d(\bx)\bigg \rangle&=\bigg\langle \mathcal{A},\mathcal{T}_{d}(\bx)\bigg \rangle\\ 
&=\bigg\langle \mathcal{A},\mathcal{T}_{d}\bigg(\sum_{i=1}^{n}\lambda_i\bv_i\bigg) \bigg\rangle \\
&=  \sum_{i_1, \ldots, i_d =1}^n {\lambda_{i_1}\lambda_{i_2} \cdots \lambda_{i_d}} \bigg\langle \mathcal{A}, \bv_{i_1}\otimes \bv_{i_2}\otimes \cdots\otimes \bv_{i_d}\bigg\rangle \ge 0
\end{align*}
As both $\sum_{i_1, \ldots, i_d =1}^n {\lambda_{i_1}\lambda_{i_2} \cdots \lambda_{i_d}} $ and $\langle \mathcal{A}, \bv_{i_1}\otimes \bv_{i_2}\otimes \cdots\otimes \bv_{i_d}\rangle$ are non-negative.
\eop \\
\noindent
For the standard simplex $\Delta^s=conv\{\be_1,\be_2,\cdots,\be_n\}$ by the implication of Lemma \ref{lemma1}, the tensor $\mathcal{A}$ is copositive if $0\le \bigg\langle \mathcal{A} \ , \mathcal{T}_{d}\big(\sum_{i=1}^{n}\lambda_i\be_i \big) \bigg\rangle=a_{i_1i_2 \cdots i_d}$, which establishes the fact that any entry wise non-negative tensor is copositive.  looking at the simplicial partition which is stated as under : \\
Let $\Delta$ be a simplex in $\Re^n$ a family $\mathcal{P}=\{\Delta^1,\Delta^2,\cdots,\Delta^m\}$ of sub-simplexes satisfying
\begin{align*}
\Delta=\bigcup_{i=1}^{m}\Delta^i ~ ~ \text{and} ~ ~ int(\Delta^i) \bigcap int(\Delta^j)=\phi ~~for~~i \ne j
\end{align*}
is said to be simplicial partition of $\Delta.$ The set of all vertices's in $\mathcal{P}$ is denoted by $V_\mathcal{P}$ and the set of all edges in $\mathcal{P}$ by $E_\mathcal{P}$. One of the convenient approach to partition any simplex is the radial subdivision of
$\Delta$; choosing an arbitrary $\bu\in\Delta\backslash \{\bv_1,\bv_2,\cdots,\bv_n\}$ such that, $\bu=\frac{\bv_{i}+\bv_{i+1}}{2}$
The sub-simplex $\Delta^i$ is obtained by replacing the vertex $\bv_i$ in $\Delta$ with $u$; that is,
\begin{align*}
\Delta^i=conv\{\bv_1,\cdots,\bv_{i-1},\bu,\bv_{i+1},\cdots,\bv_n\}
\end{align*}
For brevity of notation we define $\mathcal{T}_{d}(\bx,\by)=\bx^{\alpha} \otimes \by^{d-\alpha}$.  We describe the sufficient conditions for copositivity, in the following theorem which is generalization of [Theorem 3,\cite{Bundfuss2009}];
\begin{theorem}\label{theorem1}
	Let $\mathcal{A}$ be a $d^{th}$-order, $n$-dimensional symmetric tensor, and $\mathcal{P}=\{\Delta^1,\cdots,\Delta^m\}$ be a simplicial partition of $\Delta^s$; if
	\[ \bigg\langle \mathcal{A},{\mathcal{T}_{d}(\bu,\bv)} \bigg \rangle \ge 0 ~  \text{ for  each  pair} ~ \{\bu,\bv\}\in E_{\mathcal{P}} ~ \text{and }~ \bigg \langle \mathcal{A} , \mathcal{T}_{d}(\bv) \bigg \rangle\ge 0  ~ \forall ~ \bv\in V_\mathcal{P}\] 
	then $\mathcal{A}$ is copositive.
\end{theorem}
\noindent
\Pf
Let $\mathcal{P}=\{\Delta^1, \Delta^1, \cdots,\Delta^m\}$ be a simplicial partition of $\Delta^s$ then, for each
$\bu,\bv \in \Delta^s$ there exist some $\Delta^i,\Delta^j\in\mathcal{P}$ such that $\bu \in \Delta^i$ and $\bv \in \Delta^j$, by hypothesis it is true that, for all possible pairs of vertices's in the partition $\mathcal{P}$ of standard simplex $\Delta^s$ we have;  $\big \langle \mathcal{A},{\mathcal{T}_{d}(\bu,\bv)} \big \rangle \ge 0$. The standard simplex $\Delta^s$ is topologically equivalent to $\Re^n_+$, since there exists a one-one, onto mapping $ \phi:\Re^n_+ \rightarrow \Delta^s $ defined as:
\begin{align*}
\phi(\bx)=\frac{\sum_{i=1}^{n}{x_{i}}\be_{i}}{\|\bx\|_{1}}~~\text{where}~~x_i \in \Re_+
\end{align*}
Thus, corresponding to each $\bx \in \Re^n_+$ there exists distinctly unique  $\phi(\bx)\in\Delta^s$, such that;
\begin{align*}
 \bigg\langle\mathcal{A} \ , \ \mathcal{T}_{d} \big({\bx}\big) \bigg \rangle &=\|\bx\|_{1} \bigg\langle\mathcal{A} \ , \ \mathcal{T}_{d} \bigg(\frac{\bx}{\|\bx\|_{1}}\bigg) \bigg \rangle \\ 
 &= \|\bx\|_{1}  \bigg\langle\mathcal{A} \ \ , \ \mathcal{T}_{d} \big(\phi(\bx)\big) \ \bigg\rangle \\ 
 &=\|\bx\|_{1} \bigg\langle\mathcal{A} \ \ , \ \mathcal{T}_{d} \big(\bu \big)  \bigg\rangle  \ge  0
\end{align*} 
Thus, $\bigg\langle\mathcal{A} \ , \ \mathcal{T}_{d} \big({\bx}\big) \bigg\rangle \ge 0 ~\forall~\bx \in \Re^n_+$ implies that, $\mathcal{A}$ is copositive tensor.
\eop \\
\noindent
 we define the diameter $\delta(\mathcal{P})$ of the partition $\mathcal{P}$ as follows:
\begin{align*}
\delta(\mathcal{P}):=\underset{\{\bu,\bv\}\in E_\mathcal{P}}{\max} \|\bu-\bv\|
\end{align*}
As the diameter $\delta(\mathcal{P})$ tends to zero, the partition gets finer and finer, eventually members of strictly copositive tensor's cone are captured. For the limiting case we state the necessary condition for a tensor to be strictly copositive.
\begin{theorem}\label{theorem2}
	Let $ \mathcal{A}\in \mathcal{S}_{n,d}$ be a strictly copositive tensor, then for every finite simplicial partition $\mathcal{P}$ of $\Delta$ there exists an $\varepsilon>0$ with $\delta(\mathcal{P})\le\varepsilon$ such that;
	\begin{align*}
	\bigg\langle \mathcal{A},{\mathcal{T}_{d}(\bu, \bv)}\bigg\rangle & \ge 0  ~~\forall~~\{\bu,\bv\}\in E_{\mathcal{P}}~~ \\ 
	\bigg\langle \mathcal{A},\mathcal{T}_{d}(\bv) \bigg\rangle & \ge 0~~ \forall~~ \bv\in V_\mathcal{P}
	\end{align*}
\end{theorem}
\noindent
\Pf
Let $\mathcal{A} \in \mathcal{S}_{n,d}$ be a strictly copositive tensor, then the associated polynomial form $f_\mathcal{A}(\bx) > 0~~ \forall ~ \bx \in \Re^n_+$. Since $\Delta$ is a compact subspace of $\Re^n_+$, therefore by continuity condition it implies that, for each $\bx,\by\in\Delta$ there exists an $\varepsilon_{\bx}\ge0$ such that;
\begin{align*}
f_\mathcal{A}(\bx,\by)=\bigg\langle \mathcal{A},\mathcal{T}_{d}{(\bx,\by)} \bigg\rangle > 0  ~~\forall~~\|\bx-\by\| \le \varepsilon_{\bx}
\end{align*}
By uniform continuity of polynomials on the compact space $\Delta$, it implies that; $\varepsilon:=\underset{\bx\in\Delta}{\inf}\varepsilon_{\bx}~>0$.
For the simplicial partition $\mathcal{P}$  of simplex $\Delta$ with $\delta(\mathcal{P}) \le \varepsilon$, choose an arbitrary sub-simplex $\Delta^k$ and $\bx^{(k)},\by^{(k)}\in\Delta^k$ it implies that $\|\bx^{(k)}- \by^{(k)}\|\le\varepsilon$, which further implies that for all possible pairs of vertices's $\bx^{(k)},\by^{(k)}$  in  $\Delta$, we have;  
\begin{align*}
&f_\mathcal{A}(\bx^{(k)},\by^{(k)}) \ge 0 ~~as ~~k\rightarrow \infty\\
\implies & f_\mathcal{A}(\bx,\by) \ge 0 ~~\forall~~\{\bx,\by\} \in E_\mathcal{P}~~and~~f_\mathcal{A}(\bx) \ge 0 ~~\forall~~\bx \in V_\mathcal{P}
\end{align*}
\eop
Consequently, for any partition $\mathcal{P}$ and  by Theorems \ref{theorem1} and \ref{theorem2}, it is natural to define inner polyhedral approximations for the copositive cone $\mathcal{C}_{n,d}$ as follows:
\begin{small}
	\begin{align*}
\mathcal{I}_{n,d}^{\mathcal{P}}:=\bigg\lbrace\mathcal{A} \in {\mathcal{S}}_{n,d}: \bigg\langle \mathcal{A},{\mathcal{T}_{d}(\bu,\bv)} \bigg\rangle \ge 0  ~\forall~\{\bu,\bv\}\in E_{\mathcal{P}} ~\&~\bigg\langle \mathcal{A},\mathcal{T}_{d}(\bv) \bigg\rangle\ge 0 ~~\forall~~\bv\in V_{\mathcal{P} } \bigg\rbrace 
\end{align*}
\end{small}
\noindent
The inner cone approximation  $\mathcal{I}_{n,d}^{\mathcal{P}}$ corresponding to the partition $\mathcal{P}$, and for two simplicial partitions $\mathcal{P}_1$ and $\mathcal{P}_2$ of simplex $\Delta$. The partition $\mathcal{P}_2$ is said to be refinement of $\mathcal{P}_1$ if for each sub-simplex $\Delta^k\in\mathcal{P}_1$ there exists a
subset $\mathcal{P}_{\Delta^k}\subseteq \mathcal{P}_2$ which is simplicial partition of $\Delta^k$. In the subsequent part of this section, we discuss several properties of $\mathcal{I}_{n,d}^{\mathcal{P}}$.
\begin{lemma}\label{lemma2}
	The cone $\mathcal{I}_{n,d}^{\mathcal{P}}$ is an inner approximation for $\mathcal{C}_{n,d},$ that is, $\mathcal{I}_{n,d}^{\mathcal{P}}\subseteq \mathcal{C}_{n,d}$.
\end{lemma}
\noindent
\Pf
Let $\mathcal{A} \in \mathcal{I}_{n,d}^{\mathcal{P}}$ be an arbitrary tensor, then for each $\bx\in\Delta$ there exist a sub-simplex $\Delta^k$ such that; $\bx\in\Delta^k$ and by continuity conditions it implies that for $\delta(\mathcal{P})\le\varepsilon$ and $\bx,\by\in \Delta^k$ we have;
\begin{align*}
&f_\mathcal{A}(\bx,\by)=\bigg\langle \mathcal{A},\mathcal{T}_{d}(\bx,\by) \bigg\rangle \ge 0 ~~whenever~~\|\bx-\by\| \le \varepsilon \\
\implies & f_\mathcal{A}(\bx,\by) \ge 0 ~~\forall~~ \bx,\by\in \Re^n_+ \\
\implies & f_\mathcal{A}(\bx) \ge 0 ~~\forall~~ \bx \in \Re^n_+
\end{align*}
Hence $\mathcal{A}$ is copositive tensor, that is $\mathcal{I}_{n,d}^{\mathcal{P}}\subseteq \mathcal{C}_{n,d}$.
\eop
 we consider a sequence $\{\mathcal{A}_k\}$ of tensors in $\mathcal{I}_{n,d}^{\mathcal{P}}$, with reference to Theorem \ref{theorem2}; we see that, for any partition $\mathcal{P}_k$ as $\delta(\mathcal{P}_k)\rightarrow 0$ we have; 
\begin{align*}
&\lim\limits_{k \rightarrow \infty}\bigg\langle\mathcal{A}_k,{\mathcal{T}_{d}(\bx^{k})} \bigg\rangle \ge 0 ~~\forall~\bx^k\in\Delta^k\\
& \implies \bigg\langle\mathcal{A},\mathcal{X} \bigg \rangle \ge 0 ~~\forall~\bx \in \Delta^k \text{ since } \Delta ~~is~compact
\end{align*}
Thus, $\mathcal{A}\in \mathcal{I}_{n,d}^{\mathcal{P}}$, which implies that $\mathcal{I}_{n,d}^{\mathcal{P}}$ is closed and convex as well. Moreover, for every finite partition $\mathcal{P}$ of $\Delta$ the total number of vertices's and edges are also finite, therefore to determine that, a tensor $\mathcal{A}\in \mathcal{I}_{n,d}^{\mathcal{P}}$ it required to solve finitely many inequalities. Henceforth such tensor cones can be generated by finite subset $\mathcal{M}\subseteq \mathcal{I}_{n,d}^{\mathcal{P}}$, which establishes the fact that $\mathcal{I}_{n,d}^{\mathcal{P}}$ is polyhedral cone. In the following lemma we discuss the containment relation among inner approximations based on partition $\mathcal{P}_1$ and its refinement $\mathcal{P}_2$.
\begin{lemma}\label{lemma3}
	Let $\mathcal{P}$, $\mathcal{P}_1$, and $\mathcal{P}_2$ be simplicial partitions of $\Delta$. If $\mathcal{P}_2$ is a refinement of $\mathcal{P}_1$, then $\mathcal{I}_{n,d}^{\mathcal{P}_1} \subseteq \mathcal{I}_{n,d}^{\mathcal{P}_2}$		
\end{lemma}
\noindent
\Pf
Let $\mathcal{P}_2$ be a refinement of $\mathcal{P}_1$, then there exists a sub-simplexes $\Delta^{k_1}\in\mathcal{P}_1$ and $\Delta^{k_2}\in\mathcal{P}_2$ such that $\Delta^{k_2}\subseteq\Delta^{k_1}$.  we consider an arbitrary tensor $\mathcal{A}\in\mathcal{I}_{n,d}^{\mathcal{P}_1}$ , and $\bx,\by \in \Delta^{k_2}$. As both $\bx$ and $\by$ can be expressed as convex combination of the vertices's $\bu^{k_1}_i\in\Delta^{k_1}$:
\[	\bx=\sum_{i=1}^{n}{\lambda_i \bu^{k_1}_i}~~where~~\sum_{i=i}^{n}\lambda_i=1; ~~\lambda_i\ge0 \]
\[\by=\sum_{i=1}^{n}{\theta_i \bu^{k_1}_i}~~where~~\sum_{i=i}^{n}\theta_i=1;~~\theta_i\ge0 \]
Since for each pair $\bu^{k_1}_i,\bu^{k_1}_j$ of vertices's in $\Delta^{k_1}$ we have, $\bigg\langle \mathcal{A},\mathcal{T}_{d}{(\bu^{k_1}_i,\bu^{k_1}_j)} \bigg\rangle \ge 0$, which implies that; 
\begin{small}
\begin{align*}
f_{\mathcal{A}}(\bx,\by)&=\bigg\langle \mathcal{A} \ \ , \ \ \underset{\alpha-times}{\underbrace{\bx\otimes \cdots\otimes \bx }} \ \otimes \ \underset{({d-\alpha})-times}{\underbrace{\by\otimes \cdots\otimes \by}} \bigg\rangle
\\
&=\begin{tiny}
\underset{j_1, \ldots, j_{(d-\alpha)}=1}
{\sum_{i_1, \ldots, i_\alpha =1}^n} 
({\lambda_{i_1} \cdots \lambda_{i_\alpha}}{\theta_{j_1} \cdots \theta_{j_{(d-\alpha)}}} ) \end{tiny}
\bigg\langle \mathcal{A}, \bu_{i_1}\otimes \cdots\otimes \bu_{i_\alpha}\otimes \bu_{j_1}\otimes \cdots\otimes \bu_{j_{(d-\alpha)}} \bigg\rangle \\ 
& \ge  \ \ \ 0
\end{align*}
\end{small}
Hence, $\mathcal{A}\in\mathcal{I}_{n,d}^{\mathcal{P}_2}$, which further implies that, $\mathcal{I}_{n,d}^{\mathcal{P}_1} \subseteq \mathcal{I}_{n,d}^{\mathcal{P}_2}$.
\eop \\
\noindent
The sequence $\{\mathcal{P}_k\}$ of simplicial partition  yields a system of polyhedral inner approximations which approximates the copositive cone precisely.
\begin{theorem}\label{theorem3}
	Let the sequence $\{\mathcal{P}_k\}$ of simplicial partitions of $\Delta$ with $\delta(\mathcal{P}_k)\rightarrow0$. Then we have:
	\[int(\mathcal{C}_{n,d})\subset\underset{k\in\mathbb{N}}{\bigcup} \mathcal{I}_{\mathcal{P}_k} = \mathcal{C}_{n,d} \]
\end{theorem}
\noindent
\Pf
Suppose $\mathcal{A}\in int(\mathcal{C}_{n,d})$ be an arbitrary tensor, then $\mathcal{A}$ is strictly copositive. By
Theorem \ref{theorem2} it implies that, there exists a partition $\mathcal{P}_{k_0}$ such that $\mathcal{A}\in \mathcal{I}_{n,d}^{\mathcal{P}_{k_0}}$ where $k_{0}\in\mathbb{N}$; therefore we have, \[\mathcal{A}\in\underset{k\in\mathbb{N}}{\bigcup} \mathcal{I}_{n,d}^{\mathcal{P}_{k}}~~\implies~~ int(\mathcal{C}_{n,d})\subset\underset{k\in\mathbb{N}}{\bigcup} \mathcal{I}_{n,d}^{\mathcal{P}_{k}}\]
However, for any $\mathcal{A} \in \underset{k\in\mathbb{N}}{\bigcup} \mathcal{I}_{n,d}^{\mathcal{P}_{k}}$ there exists some partition $\mathcal{P}_{k}$ such that \[\bigg\langle \mathcal{A}, \mathcal{T}_{d}(\bv) \bigg\rangle = 0 ~~\text{for}~~\bv\in V_{\mathcal{P}_{k}} ~\implies~\mathcal{A} \notin int(\mathcal{C}_{n,d}) \] Moreover, by Lemma \ref{lemma2} we have, $\mathcal{I}_{n,d}^{\mathcal{P}_{k}}\subset\mathcal{C}_{n,d}~~\forall~~k\in\mathbb{N}$ which implies $\underset{k\in\mathbb{N}}{\bigcup} \mathcal{I}_{n,d}^{\mathcal{P}_{k}}\subseteq\mathcal{C}_{n,d}$.\\
For an arbitrary $\mathcal{A} \in \mathcal{C}_{n,d}$ we have;
\begin{align*}
&\bigg\langle \mathcal{A}, \mathcal{T}_{d}(\bx) \bigg\rangle \ge 0 ~~\text{for}~~\bx\in \Re^{n}_{+} \\ 
\implies & \|\bx\|_{1}\bigg\langle \mathcal{A} , \mathcal{T}_{d}\bigg(\frac{\bx}{\|\bx\|_{1}}\bigg) \bigg\rangle \ge 0 ~~\text{for}~~{\frac{\bx}{{\|\bx\|_{1}}}} \in \Delta  \\
\end{align*}
hence, there exists a partition $\mathcal{P}_k$ of simplex  $\Delta$ such that;
\[f_{\mathcal{A}} \bigg(\frac{\bx}{{\|\bx\|_{1}}}\bigg) \ge 0  \ \ \text{for all} \ \frac{\bx}{{\|\bx\|_{1}}} \in V_{\mathcal{P}_k} ~\text{ and } f_{\mathcal{A}} \bigg(\frac{\bx}{{\|\bx\|_{1}}},\frac{\by}{{\|\by\|_{1}}}\bigg) \ge 0 \text{  for all  } \bigg\lbrace \frac{\bx}{{\|\bx\|_{1}}},\frac{\by}{{\|\by\|_{1}}} \bigg\rbrace \in E_{\mathcal{P}_k} \] 
Thus, $\mathcal{A} \in \mathcal{I}_{n,d}^{\mathcal{P}_k}$ as well, therefore we have the result
\begin{align*}
\underset{k\in\mathbb{N}}{\bigcup} \mathcal{I}_{n,d}^{\mathcal{P}_{k}} = \mathcal{C}_{n,d}
\end{align*}
\eop
\subsection{Containment relations among $\mathcal{C}_{n,d}^{(r)}$ , $\mathcal{K}_{n,d}^{(r)}$ and $\mathcal{I}_{n,d}^{\mathcal{P}_r}$}
We present the inclusion relation between the cones $\mathcal{C}_{n,d}^{(r)}$, $\mathcal{K}_{n,d}^{(r)}$ and $\mathcal{I}_{n,d}^{\mathcal{P}_r}$ as a proposition:
\begin{proposition} Let $\mathcal{P}_r$ be a simplicial partition of the simplex $\Delta$, then the $r^{th}$ level inner approximation hierarchies $\mathcal{C}_{n,d}^{(r)}$ ,  $\mathcal{K}_{n,d}^{(r)}$ and $\mathcal{I}_{n,d}^{\mathcal{P}_r}$  for the copositive tensor cone $\mathcal{C}_{n,d}$ has the following inclusion relations:  
	\begin{align*}
	\mathcal{C}_{n,d}^{(r)} \subseteq \mathcal{K}_{n,d}^{(r)}~~and~~\mathcal{C}_{n,d}^{(r)} \subseteq \mathcal{I}_{n,d}^{\mathcal{P}_r} ~~for ~~all~~r \in \{0,1,2,\cdots\} 
	\end{align*}
\end{proposition} 
\noindent
\Pf Let $ \mathcal{A} \in \mathcal{C}_{n,d}^{(r)}$ be an arbitrary tensor, then the associated polynomial $P^{(r)}(\by)$ in (\ref{multicoeff}) allows a sum-of-square decomposition, thus  $ \mathcal{A} \in \mathcal{K}_{n,d}^{(r)}$, henceforth, $\mathcal{C}_{n,d}^{(r)} \subseteq \mathcal{K}_{n,d}^{(r)}~\forall~r$. \\
To show that the tensor $ \mathcal{A}$ belongs to $r^{th}$ level cone $\mathcal{I}_{n,d}^{\mathcal{P}_r}$ as well, we take an arbitrary $\bv \in V_{\mathcal{P}_r}$, which implies that $(r+d)\bv \in \mathbb{I}^n(r+d)$. Therefore for any pair of vertices's $\bu,\bv \in V_{\mathcal{P}_r}$ we have; 
\begin{align*}
\bigg\langle \mathcal{A} \ , \ \mathcal{T}_{d}(\bu,\bv) \bigg\rangle 
&=\bigg\langle \mathcal{A},\mathcal{T}_{d}\bigg(\frac{\bx}{r+d},\frac{\by}{r+d}\bigg) \bigg\rangle \\
&=\frac{1}{(r+d)^d} \bigg\langle \mathcal{A},\mathcal{T}_{d}(\bx,\by) \bigg\rangle \\
&=\frac{1}{(r+d)^d} \bigg(\bigg\langle \mathcal{A},\mathcal{T}_{d}(\bx) \bigg\rangle + \bigg\langle \mathcal{A},\mathcal{T}_{d}(\by) \bigg\rangle- \bigg\langle \mathcal{A},\mathcal{T}_{d}(\bx-\by) \bigg\rangle \bigg)\\
& \ge \frac{1}{(r+d)^d} \bigg(\bigg\langle \mathcal{A},Diag(\bx) \bigg\rangle + \bigg\langle \mathcal{A},Diag(\by) \bigg\rangle- \bigg\langle \mathcal{A},Diag(\bx-\by) \bigg\rangle\bigg)\\
&\ge \frac{1}{(r+d)^d} \sum_{i=1}^{n} \bigg(a_{(i)^d} x_i + a_{(i)^d} y_i- a_{(i)^d} (x_i-y_i) \bigg) \\
&\ge \frac{1}{(r+d)^d} \sum_{i=1}^{n} \big(2a_{(i)^d} y_i\big) \ge 0
\end{align*}
Thus $\bigg\langle \mathcal{A}  ,  \mathcal{T}_{d}(\bu,\bv) \bigg\rangle \ge 0~~\forall~~ \bu,\bv \in V_{\mathcal{P}_r}$, which further implies $\bigg\langle \mathcal{A} , \mathcal{T}_{d}(\bv) \bigg\rangle  \ge 0 ~~\forall~~ \bv \in V_{\mathcal{P}_r}$. Hence, $\mathcal{A} \in \mathcal{I}_{n,d}^{\mathcal{P}_r}$ for all $r$.
\eop

Note that, neither $\mathcal{I}_{n,d}^{\mathcal{P}_r} \subseteq \mathcal{K}_{n,d}^{(r)}$ nor $\mathcal{K}_{n,d}^{(r)} \subseteq \mathcal{I}_{n,d}^{\mathcal{P}_r}$. 
\begin{example}
	Let $\mathcal{A}\in \mathcal{K}_{3,6}^{(0)}$ be any tensor with entries given as follows: 
	\begin{align*}
	a_{i_1i_2i_3i_4i_5i_6}= \begin{cases}
		\begin{matrix}
		\ 1 & \ if \ i_{j}=i_{k}~\forall~j,k \in \{1,2,3,4,5,6\} \\
		 \ 2 & \ if \ i_{j}=1~\forall~j\in\{1,2,3\} ~and~ i_{K}=2~\forall~k \in \{4,5,6\} \\
		\ 2 & \ if \ i_{j}=1~\forall~j \in \{1,2,3\} ~and~ i_{K}=3~\forall~  k \in \{4,5,6\}  \\
		-2 & \ if \ i_{j}=2~\forall~j \in \{1,2,3\} ~and~ i_{K}=3~\forall~ k \in  \{4,5,6\} \\ 
		0 & \ otherwise
		\end{matrix}	
	\end{cases}
	\end{align*}
	then the associated polynomial form,  $f_{\mathcal{A}}(x,y,z)=x^6+y^6+z^6+2x^3y^3+2x^3z^3-2y^3z^3$  must allow an sos decomposition. , since not all elements of $\mathcal{A}$ are non-negative, therefore $ \mathcal{A} \notin \mathcal{I}_{3,6}^{\mathcal{P}_0}$.
\end{example}
\begin{example}
	Let $\mathcal{A} \in \mathcal{I}_{2,2}^{\mathcal{P}_1}$ be an arbitrary tensor, then the following conditions on elements of $\mathcal{A}$ are true:
	\begin{align*}
	&a_{ii} \ge 0 ~~for~i\in \{1,2\}\\
	&a_{11}+a_{12} \ge 0\\
	&a_{11}+a_{12} \ge 0
	\end{align*}
	the associated polynomial of $\mathcal{A}$ is $f_\mathcal{A}(\bx)=a_{11}x_1^2+a_{22}x_2^2+2a_{12}x_1x_2 ~~for~~\bx\in\Re^n_+$, by substituting ${\bx}={\by} \circ {\by}~~for~\by\in\Re^n$, which further leads to the polynomial 	
	\begin{align*}
	P^{(1)}(\by)&=(a_{11} y_{1}^4+a_{22} y_{2}^4+2a_{12}y_{1}^2 y_{2}^2) (y_{1}^2+y_{2}^2)\\
	&=\begin{pmatrix}
	y_1^2\\
	y_2^2\\
	y_1y_2
	\end{pmatrix}  ^{T} 
	\begin{pmatrix}
	a_{11}y_1^2 & 0 & a_{12}y_1y_2\\
	0 & a_{22}y_2^2 & a_{12}y_1y_2\\
	a_{12}y_1y_2 & a_{12}y_1y_2 & a_{11}y_1^2+a_{22}y_2^2
	\end{pmatrix} 
	\begin{pmatrix}
	y_1^2\\
	y_2^2\\
	y_1y_2
	\end{pmatrix}  \\
	&=V^{T}Q(y_1,y_2)V 
	\end{align*}
	 since, $Q(y_1,y_2)$ is not a PSD matrix; for if taking $y_1=y_2=1$ and for the vector $z^T=(-10,-10,1)$, we have; $z^T Q(1,1)z=101(a_{11}+a_{22})-40a_{12}$ which is negative for $a_{11}=a_{22}=\frac{1}{101}~~and~~a_{12}=1$, that is $z^T Q(1,1)z=-38$. Henceforth, $P^{(1)}(\by)$ don't allow sum-of-square decomposition. Thus $\mathcal{A}\notin \mathcal{K}_{2,2}^{(1)}$.
\end{example}
\subsection{Outer Approximations for Copositive Cone}
In this section we discuss, the outer approximations $ \mathcal{O}_{n,d}^\mathcal{P}$ and $\mathcal{O}_{n,d}^{(r)}$ for copositive cone of tensors. These approximations contains copositive cone $\mathcal{C}_{n,d}$.  
\subsubsection{Outer Approximations based on Simplicial Partition}
For any simplicial partition $\mathcal{P}$  of $\Delta$ the cone $\mathcal{O}_{n,d}^\mathcal{P}$ is stated as follows;
\begin{align*}
\mathcal{O}_{n,d}^\mathcal{P}=\bigg\lbrace\mathcal{A} \in \mathcal{S}_{n,d}: \big\langle \mathcal{A},\mathcal{T}_{d}(\bv) \big\rangle \ge 0 ~~\forall~\bv\in V_{\mathcal{P}} \bigg\rbrace
\end{align*} 
approximates the copositive cone from outside. Since for each tensor $\mathcal{A} \in \mathcal{C}_{n,d}$; the associated polynomial form $f_\mathcal{A}(\bx)=\big\langle \mathcal{A},\mathcal{T}_{d}(\bx) \big\rangle \ge~0~\forall~\bx\in\Re^n_+$. For the collection of all vertices's $V_\mathcal{P}\subseteq\Re^n_+$, we have, $\big \langle \mathcal{A}, \mathcal{T}_{d}(\bx) \big \rangle \ge~0~\forall~\bx\in V_{\mathcal{P}}\subseteq\Re^n_+$. Consequently, the tensor $\mathcal{A}\in \mathcal{O}_{n,d}^\mathcal{P}$, which implies $\mathcal{C}_{n,d} \subseteq \mathcal{O}_{n,d}^{\mathcal{P}}$. Moreover, if the partition $\mathcal{P}=\{\Delta^{s}\}$ then  $\mathcal{O}_{n,d}^{\{\Delta^{s}\}}$ carries all those tensors whose diagonal entries are non-negative. It is well known fact that the diagonal entries of copositive tensors are necessarily non-negative.\\
Clearly, the cone $\mathcal{O}_{n,d}^{\mathcal{P}}$ is closed and convex, and it possesses a polyhedral geometry, since for any finite partition $\mathcal{P}$ of simplex $\Delta$ the total number of vertices's $V_{\mathcal{P}}$ are finite, therefore to find an element of such cone requires to solve finitely many inequalities. , if $\mathcal{P}, \mathcal{P}_1,$ and $\mathcal{P}_2$ are the simplicial partitions of $\Delta$ then the inclusion relation between the outer approximations based on these partitions is presented as a lemma, whose proof is analogous to [Lemma 8,\cite{Bundfuss2009}]; for the sake of completeness we give the proof.
\begin{lemma}\label{lemma3'}
	For simplicial partitions $\mathcal{P}$, $\mathcal{P}_1$, and $\mathcal{P}_2$ of $\Delta$, if $\mathcal{P}_2$ is a refinement of $\mathcal{P}_1$, then $\mathcal{O}_{n,d}^{\mathcal{P}_2} \subseteq \mathcal{O}_{n,d}^{\mathcal{P}_1}$.		
\end{lemma}
\noindent
\Pf
Let $\mathcal{P}_2$ be a refinement of $\mathcal{P}_1$; so there exists a sub-simplex $\Delta^{k_1}\in\mathcal{P}_1$ and another sub-simplex $\Delta^{k_2}\in\mathcal{P}_2$ such that $\Delta^{k_2} \subseteq \Delta^{k_1} $, which implies that, $ V_{\mathcal{P}_1}\subseteq V_{\mathcal{P}_2}$, therefore the set of all inequalities defining $\mathcal{O}_{n,d}^{\mathcal{P}_1}$ is a subset of the set of inequalities defining $\mathcal{O}_{n,d}^{\mathcal{P}_2}$; thus, $\mathcal{O}_{n,d}^{\mathcal{P}_2} \subseteq \mathcal{O}_{n,d}^{\mathcal{P}_1}$.
\eop
The sequence $\{\mathcal{O}_{\mathcal{P}_k}\}$ of outer approximations  converges to $\mathcal{C}_{n,d}$ as the radius $\delta(\mathcal{P}_k)\rightarrow0$, as stated in the following theorem.
\begin{theorem}\label{theorem4}
	For the sequence $\{\mathcal{P}_k\}$ of simplicial partitions of $\Delta$ with $\delta(\mathcal{P}_k)\rightarrow 0$; we have:
	\begin{align*}
	\mathcal{C}_{n,d}=\underset{k\in\mathbb{N}}{\bigcap} \mathcal{O}_{n,d}^{\mathcal{P}_k} 
	\end{align*}	
\end{theorem}
\noindent
\Pf
Since we know that; $\mathcal{C}_{n,d} \subseteq \mathcal{O}_{n,d}^{\mathcal{P}_k}~~\forall k\in\mathbb{N}$ therefore, $\mathcal{C}_{n,d} \subseteq \underset{k\in\mathbb{N}}{\bigcap}\mathcal{O}_{n,d}^{\mathcal{P}_k}$. 
To establish the reverse inclusion, we assume on the contrary that, $\mathcal{A}\notin\mathcal{C}_{n,d}$ which implies that, for some $\bv\in\Delta$ we have $\big\langle \mathcal{A}, \mathcal{T}_{d}(\bv) \big\rangle < 0 $, so by continuity property, there exists an $\varepsilon$-neighborhood $N_\varepsilon(\bv)$ such that; 
\[\big\langle \mathcal{A}, \mathcal{T}_{d}(\bv) \big\rangle<0 ~~\forall~\bv\in N_\varepsilon(\bv)\]
 if the diameter $\delta(\mathcal{P})<\varepsilon$, then there exist a simplex $\Delta^{k}\in\mathcal{P}$ and $\bu^{k}\in\Delta^{k}$such that; $\|\bu-\bv\|<\varepsilon$ which implies that; $\bu^{k}\in N_\varepsilon(\bv)$, which further implies;
\begin{align*}
\big \langle \mathcal{A}, \mathcal{T}_{d}(\bu^k) \big \rangle <0 ~\forall~\bu^{k}\in N_{\varepsilon}(\bv)\subseteq \Delta^k ~\text{which is a contradiction therefore}~ \ \mathcal{A}\notin \mathcal{O}_{n,d}^{\mathcal{P}_k}
\end{align*}
Henceforth, $\mathcal{A}\notin \underset{k\in\mathbb{N}}{\cap}\mathcal{O}_{n,d}^{\mathcal{P}_k}$, which implies that, $\mathcal{C}_{n,d} \supseteq \underset{k\in\mathbb{N}}{\cap}\mathcal{O}_{n,d}^{\mathcal{P}_k}$.
\eop
\subsection{Rational Griding based Outer Approximations}
The regular grid $\Delta_{n}^{(r)}$ of rational points on the simplex $\Delta$ for each $r\in \{0,1,2,\cdots\}$ is stated as follows:  
\[ \Delta_{n}^{(r)}:=\lbrace \bx\in\Delta:(r+2)\bx\in\mathbb{N}^n_0 \rbrace \] 
 For each $r$ the grid $\Delta_{n}^{(r)}$  provides a finite discretization of the simplex $\Delta$ consisting of rational points. The cardinality of each grid $\Delta_{n}^{(r)}$ is given as follows: 
\[ |\Delta_{n}^{(r)}|={n+r-1 \choose r}\]
For each $r\in \{0,1,2,\cdots\}$, let us define \[\delta_n^{(r)}:=\bigcup_{k=1}^{r}{\Delta_{n}^{(k)}} \]
 by using the above mentioned discretization, we define  another hierarchy of outer polyhedral approximations  $ \mathcal{O}_{n,d}^{(r)}$ for copositive cone $\mathcal{C}_{n,d}$ as follows: 
\begin{align*}
\mathcal{O}_{n,d}^{(r)}:=\big\lbrace\mathcal{A} \in \mathcal{S}_{n,d}:\big \langle \mathcal{A}, \mathcal{T}_{d}(\bv) \big \rangle\ge 0 ~~\forall~~\bv\in \delta_n^{(r)} \big\rbrace 
\end{align*}
Clearly, each cone $\mathcal{O}_{n,d}^{(r)}$ is a proper cone. In the following theorem, we prove that the hierarchy of outer polyhedral approximations $\mathcal{O}_{n,d}^{(r)}$ converges to the cone of copositive tensors.
\begin{theorem}\label{theorem5}
	The hierarchy of outer polyhedral approximations $\mathcal{O}_{n,d}^{(r)}$ contains the copositive cone $\mathcal{C}_{n,d}$ for each $r\in \{0,1,2,\cdots\}$, that is, $\mathcal{O}_{n,d}^0\supseteq \mathcal{O}_{n,d}^1\supseteq \cdots \supseteq \mathcal{C}_{n,d}$; with, \[\mathcal{C}_{n,d}=\bigcap_{r\in\mathbb{N}} \mathcal{O}_{n,d}^{(r)}\] 	
\end{theorem}
\noindent
\Pf
Let $\mathcal{A} \in \mathcal{C}_{n,d}$ be an arbitrary tensor then the associated polynomial form $f_\mathcal{A}(\bx) \ge 0~~\forall~\bx \in \Re^n_+$, and since $\delta_n^{(r)} \subseteq \Re^n_+$ therefore $f_\mathcal{A}(\bx) \ge 0~~\forall~\bx\in \delta_n^{(r)}$, which implies that $\mathcal{C}_{n,d}\subseteq\mathcal{O}_{n,d}^{(r)}$, which further implies that, 
\begin{align*}
\mathcal{C}_{n,d} \subseteq \bigcap_{r\in\mathbb{N}} \mathcal{O}_{n,d}^{(r)}
\end{align*}
For the inclusion $\bigcap_{r\in\mathbb{N}} \mathcal{O}_{n,d}^{(r)} \subseteq \mathcal{C}_{n,d}$, let us consider a tensor $ \mathcal{A} \in \mathcal{S}_{n,d}$ be such that $ \mathcal{A} \notin \mathcal{C}_{n,d}$, therefore for some $r$ there exist $\bu^{(r)} \in \Delta_{n}^{(r)}$ such that,  $\big\langle \mathcal{A},\mathcal{T}_{d}(\bu^{(r)}) \big\rangle < 0$. Now, perturbing those components of $\bu^{(r)}$ which are zero by a small positive value $\theta>0$ in such a way that $\bu^{(r)}>0$. By continuity property of polynomials, we know that there exists an $\varepsilon_r$ neighborhood $N_{\varepsilon_r}(\bu^{(r)})$ such that $\big\langle \mathcal{A},\mathcal{T}_{d}(\bv)\big \rangle < 0$ for all $\bv\in N_{\varepsilon_r}(\bu)$. Letting, $\varepsilon_m=min\{\varepsilon_r,min_{i=1,2,\cdots,n}u_i \}$. As the set of rationals is dense in the set of reals, so there exists a vector $\mathbf{w} \in \mathbb{Q}^n_+$, with $\|\bu^{(r)}-\mathbf{w}\|<\varepsilon_m$, which implies that $\mathbf{w}>0.$ Hence there exists some positive integer $k$ such that $\bx=\frac{\mathbf{w}}{\|\mathbf{w}\|_1} \in\delta_n^{(r)}$ for all $r \ge m$, because  $\big \langle \mathcal{A},\mathcal{T}_{d}(\bx)\big \rangle < 0$, thus  $\mathcal{A}\notin \mathcal{O}_{n,d}^{(r)} ~~\forall~ r\ge m$ which further implies that, $\mathcal{A}\notin \underset{r\in\mathbb{N}}{\bigcap} \mathcal{O}_{n,d}^{(r)}$.
\eop \\
Therefore, by using the classical partitioning of the simplex $\Delta$ through bisection of longest edge, the collection of all the vertices's in the partition $\mathcal{P}_r$ of the simplex $\Delta$ is always contained in $\delta_n^{(r)}$ for each $r$. We present the containment relation between  outer approximations in the following proposition. 
\begin{proposition}
	For simplicial partition $\mathcal{P}_r$ through bisection along the longest edge and rational discretization  $\delta_n^{(r)}$ of the simplex $\Delta$, we have $\mathcal{O}_{n,d}^{\mathcal{P}_r} \subseteq \mathcal{O}_{n,d}^{(r)}~~for~ each~ r \in \{0,1,2,\cdots\}$ 
\end{proposition}
\noindent
\Pf
Let $\mathcal{P}_r$ be the classical partitioning of the simplex $\Delta$ through bisection along the  longest edge, and $V_{\mathcal{P}_r}$ be the collection of all the vertices's in $\mathcal{P}_r$, then for each $r$ we immediately have, $V_{\mathcal{P}_r} \subseteq \delta_n^{(r)}$, which implies that
\begin{align*}
\mathcal{O}_{n,d}^{\mathcal{P}_r} \subseteq \mathcal{O}_{n,d}^{(r)} ~~\forall ~ r \in \{0,1,2,\cdots \}
\end{align*} 
\eop
\section{Conclusion}
In this article several properties for Copositive Tensor cones are proven. Moreover, a necessary and sufficient condition under which the cones $\cC_{n,d}$ and $\cS^{+}_{n,d}$ coincides has been established.  The calculation of the coefficients of higher degree polynomial is critical in the analysis of approximation hierarchies based on polynomial conditions. In this regard, procedure/formula to find the coefficients of the polynomial form is obtained, along with the representation of the inner approximation cone $\cC^{r}_{n,d}$ for copositive cone $\cC_{n,d}$. More importantly, several inner and outer approximation hierarchies along with their containment relations are also given. In future we worked towards utilizing these hierarchies for approximating polynomial optimization. Especially to recover approximation results for polynomial optimization over the simplex as obtained by De Klerk and co-authors \cite{DEKLERK2006210}, \cite{Monique2013}, \cite{DeKlerkError}, \cite{Klerk2017OnTC}.

\footnotesize
\bibliography{Present}
\bibliographystyle{SIAM}
\end{document}